\documentclass[12pt]{elsart}
\usepackage{amsmath, amsfonts, amssymb}

\let\BFseries\bfseries\def\bfseries{\BFseries\mathversion{bold}} 
\newcommand{\ind}{1\hspace{-0.098cm}\mathrm{l}}
\newcommand{\eps}{\varepsilon}
\theoremstyle{remark}
\newtheorem{re}[thm]{Remark}
\newtheorem{exa}[thm]{Example}
\newenvironment{proof}[1][] {\noindent {\bf Proof#1:} }{\hspace*{\fill}$\square$\smallskip\par}

\def\P{{\bf {\mathbb{P}}}}
\newcommand{\pr}[1]{\P\left(#1\right)}
\def\ep{{\varepsilon}}
\def\E{\mathbb{E}}\def\DD{{\bf D}} \def\EE{{\mathcal E}}
\newcommand{\norm}[1]{\left\|#1\right\|}

\def\R{\mathbb{R}}\def\N{\mathbb{N}}\def\CC{\mathcal{C}} \def\d{\mathrm{d}}
\newcommand{\cD}{\mathcal{D}}
\newcommand{\conda}{condition~(*)}
\newcommand{\dist}{\operatorname*{dist}}
\newcommand{\deq}{\stackrel{d}{=}}
\newcommand{\suptwo}[2]{\sup_{\substack{#1 \\ #2}}} 

\begin{document}
\begin{frontmatter}
\title{High resolution quantization and entropy coding of jump processes}
\author{Frank\ Aurzada},
\ead{aurzada@math.tu-berlin.de}
\author{Steffen Dereich},
\ead{dereich@math.tu-berlin.de}
\author{Michael Scheutzow}, and
\ead{ms@math.tu-berlin.de}
\author{Christian Vormoor}
\ead{christian.vormoor@devk.de}
\address{Technische Universit\"{a}t Berlin, Institut f\"{u}r Mathematik, Sekr.\ MA 7-5, Stra\ss e des 17.\ Juni 136, 10623 Berlin, Germany.}

\begin{abstract}
We study the quantization problem for certain types of jump processes. The probabilities for the number of jumps are assumed to be bounded by Poisson weights. Otherwise, jump positions and increments can be rather generally distributed and correlated. We show in particular that in many cases entropy coding error and quantization error have distinct rates. Finally, we investigate the quantization problem for the special case of $\R^d$-valued compound Poisson processes.
\end{abstract}

\begin{keyword} High resolution quantization; entropy coding; complexity; jump process; compound Poisson process; L\'{e}vy process: metric entropy

\end{keyword}
\end{frontmatter}

\section{Introduction and results}
\subsection{Statement of the problem} In this article, we study the quantization- and entropy coding problem for certain types of jump processes. Given a random variable $X$, the aim is to find a good approximation $\hat X$ to $X$ that satisfies a particular \emph{complexity constraint}.

Let $s>0$, $X$ be a random variable in a measurable space $(E,\EE)$, $\rho$ a distortion measure on $E$ (i.e.\ a measurable, symmetric function $\rho : E\times E \to \R_{\geq 0}$ with $\rho(x,y)=0$ iff $x=y$), and $r\geq 0$. Then we define the \emph{quantization error} as follows:
$$D^{(q)}( r ~|~ X,\rho,s) := \inf\left\lbrace \left( \E \min_{a\in \CC} \rho(X,a)^s\right)^{1/s}\quad : \quad \log \# \CC\leq r\right\rbrace.$$
The number $D^{(q)}$ represents the best-achievable average error  when  encoding the signal $X$ with $r$ nats. The term  `nats' is used instead of `bits', since we calculate the amount of information using the natural logarithm. Further, we investigate the \emph{entropy coding error}, which can be understood as the average error when encoding the signal $X$ using -- on average --  $r$ nats:
$$D^{(e)}( r ~|~ X,\rho,s) := \inf\left\lbrace \left( \E \rho(X,\hat{X})^s\right)^{1/s} ~:~\text{$\hat{X}$ random var.\ with $H(\hat{X})\leq r$}\right\rbrace,$$
where $H$ is the (discrete) entropy of a random variable:
$$H(X):= \begin{cases} - \sum_{x} \pr{X=x} \log \pr{X=x}& \text{$X$ discrete,} \\ \infty & \text{otherwise.} \end{cases}$$
In slight misuse of notation we also write $D^{(q)}( r ~|~ X,\norm{.},s)$ if $\rho(x,y)=\norm{x-y}$ for a norm distortion $\norm{.}$. Analogously, we deal with the entropy coding error.
We recall that $D^{(e)}( r ~|~ X,\rho,s) \leq D^{(q)}( r ~|~ X,\rho,s)$.

The problems described above arise naturally in coding theory, where for instance, the complexity of a signal has to be reduced due to capacity restrictions of a channel or simply (lossy) data compression is considered (see for instance \cite{CoTho91} for a general account on coding theory and \cite{Kol68} for a historic outline of the information constraints). Beyond these information-theoretic applications, the quantization error is tightly related to certain quadrature problems:  the quantization error can be defined equivalently as the worst-case error of a particular quadrature problem. Moreover, further quadrature problems are linked to the quantization problem via estimates involving both quantities. Recent results in that direction can be found in  \cite{CDMR08}  (see also \cite{PaPri05} for earlier results).

The analysis of the quantization- and entropy coding error started in the 40s of the 20th century. At that time research was mainly focused on finite-dimensional signals; and the numerous publications mainly appeared in the engineering literature. A mathematical account of the results for finite-di\-men\-sional signals is provided by \cite{grafluschgy}. Since about 2000 researchers are attracted by the problem in the case where the original signal is infinite-dimensional. A series of articles followed on (infinite-dimensional) random vectors $X$ that are Gaussian (see for instance \cite{DFMS03}, \cite{LuPa04}, \cite{dereichscheutzow}), diffusions (\cite{LuPa04b}, \cite{Der06a}, \cite{Der06b}), and L\'evy processes (\cite{luschgypages06}, \cite{aurzadadereich}).

In this article, we provide asymptotic estimates for the quantization- and entropy coding error for certain jump processes. The results are shown to be sharp in several cases. In contrast to the (infinite-dimensional) settings studied before, there is a qualitative difference in the (best-achievable) approximation error induced by the two constraints.

\subsection{Some notation and the model}
Let us now introduce the jump processes that we investigate in this article.

We define the space $\DD([0,1[,E)$ to be the space of all functions $f:[0,1[\to E$ that are piecewise constant and possess a finite number of jumps, where if $f$ has a jump at $t$ from the value $a\in E$ to $b\in E$, then $f(t)=b$. We endow $\DD([0,1[,E)$ with the $\sigma$-field induced by the projections. 

In the sequel, $X=(X(t))_{t\in[0,1[}$ denotes a $\DD([0,1[,E)$-valued random vector. We denote by $N_X$ the random number of jumps of $X$, let $0<Y_1<\dots<Y_{N_X}<0$ be the \emph{jump positions} of $X$, and set $Y_0=0$ and $Y_{N_X+1}=1$. Moreover, we denote by
$$
Z_i:=\rho(X(Y_{i-1}),X(Y_i))
$$
the \emph{moduli of the increments} and, in the case where $E$ is a linear space, we denote by
$$
Z^{(i)}:=X(Y_i)-X(Y_{i-1})
$$
the \emph{increments}.

As distortion measure on $\DD([0,1[,E)$ we consider \begin{equation} \rho_\DD(f,g) := \int_0^1 \rho(f(t),g(t))\, \d t,\qquad f,g \in \DD([0,1[,E), \label{eqn:distm}\end{equation} where $\rho$ is a distortion measure on $E$. It is straightforward to extend the results of this paper to the distortion measure $\rho_\DD^p(f,g) = ( \int_0^1 \rho(f(t),g(t))^p \, \d t)^{1/p}$, with $1\leq p<\infty$.

Our lower bounds require that the jump positions constitute a \emph{Poisson point process} with intensity $\lambda$. The upper bounds on the complexity are proven under weaker assumptions on $X$. Here, we only assume that the total number of jumps can be estimated against the probability weights of a Poisson random variable:
\begin{equation} \pr{N_X = k} \leq \frac{\lambda^k}{k!}\, e^{-\lambda} K,\qquad k\geq 0, \label{decrcond}\end{equation}
where $\lambda>0$ and $K\ge1$ are some fixed parameters. In particular, one can choose $K=1$, if the jump positions are induced by a Poisson point process.

Sometimes we shall also impose the following condition: 
\begin{enumerate}
\item[(*)] The jump positions are independent of the jump destinations, which means that, given the event $\{N_X=k\}$, the vector $(Y_1, \ldots, Y_k)$ is independent of the vector $(X(Y_0), \ldots, X(Y_k))$.
\end{enumerate}

Let us introduce some more notation. Firstly, we make use of the concept of metric entropy. If $\rho$ is a distortion measure on $E$ we define its covering numbers by $$N(E,\rho,\eps) :=  \min\lbrace n \in \N \,:\, \exists x_1, \ldots, x_n \in E~\forall x\in E~ \exists i~:~ \rho(x,x_i)\leq \eps\rbrace.$$
A set $\{ x_1, \ldots, x_n \}$ for which the defining property of $N$ holds is called an $\eps$-net of $E$. Note that in general one has to assume that $N(E,\rho,\eps)$ is well-defined, i.e.\ that for all $\eps>0$ there is an $\eps$-net of $E$. This is ensured if, for example, $(E,\rho)$ is a precompact metric space.
We also introduce the inverse concept of $D^{(q)}$, which we call $d^{(q)}$, given by $$d^{(q)}(\eps\,|\,X,\rho,s) := \inf\left\lbrace n\geq 1, n\in\N \,:\, D^{(q)}( \log n ~|~ X,\rho, s) \leq \eps\right\rbrace.$$
In other words, $d^{(q)}$ is the number of points needed to quantize with error at most $\eps$, i.e.\ roughly it is the inverse function of $D^{(q)}(\log (.))$.

We shall also need the notation of strong and weak asymptotics. Namely, we write $f\lesssim g$, if $\limsup f/g \leq 1$. Analogously, $f\gtrsim g$ is defined. Furthermore, $f\sim g$ means $\lim f/g = 1$. We also use $f\approx g$ if $0<\liminf f/g \leq \limsup f/g < \infty$. Finally, throughout the article $\lambda_d$ denotes the $d$-dimensional Lebesgue measure.

The paper is organized as follows. In the rest of this section we state the main results.
In Section~\ref{sec:iupper}, we state the upper bounds for both quantities under various additional assumptions. In Section~\ref{sec:ilower}, the upper bounds are complemented by corresponding lower bounds. In particular, we obtain that the upper and lower bounds are tight in many cases. Finally, Section~\ref{sec:cpp} is devoted to the particular setting where $X$ is a compound Poisson process. The proofs for the upper bounds can be found in Section~\ref{sec:pg}. There, explicit coding strategies are constructed. The lower bounds are proven in Sections \ref{sec:lb} and \ref{sec:ece} for quantization- and entropy coding, respectively.  The proofs for the lower bounds of the quantization error rely on a small ball argument, whereas the lower bounds for the entropy coding error are derived using the Shannon lower bound for a related problem.

\subsection{Upper bounds} \label{sec:iupper} Our first result concerns the case where the space $E$ has finite covering numbers. In the case where $(E,\rho)$ is a metric space, this corresponds to the assumption that $E$ is precompact.

\begin{thm} Assume that $w:=\sup_{x,y\in E} \rho(x,y)<\infty$ and that the upper box dimension
\begin{equation} \gamma:=\limsup_{\ep\to 0} \frac{\log N(E,\rho,\eps)}{\log 1/\eps}  \label{eqn:assbox}\end{equation}
is finite. Then \begin{equation} \label{eqn:q} - \log D^{(q)}( r ~|~ X,\rho_\DD,s) \gtrsim \sqrt{\frac{2}{s (1+\gamma)}\, r \log r}. \end{equation} \label{thm:q}\end{thm}

\begin{thm} Assume that $N(E,\rho,\eps)<\infty$ for all $\eps>0$ and that $w:=\sup_{x,y\in E} \rho(x,y)<\infty$.
\begin{enumerate}
\item[(a)] For all $r>r_0=r_0(\lambda)$,
\begin{equation} \label{eqn:e} D^{(e)}\left(  K \left( \lambda r + (\lambda+1) \log N(E,\rho,e^{-r}) \right) ~|~ X,\rho_\DD,s\right) \leq C_{s} (w+1) K^{1/s} e^{- r},\end{equation} where the constant $C_s$ depends on $s$ only, and $K$ and $\lambda$ are the constants from (\ref{decrcond}).
\item[(b)] In particular, if the jump positions are distributed according to a Poisson point process with rate $\lambda$, we obtain for $r>r_0=r_0(\lambda)$ \begin{equation} \label{eqn:epoisson} D^{(e)}(   \lambda r + (\lambda+1) \log N(E,\rho,e^{-r}) ~|~ X,\rho_\DD,s) \leq C_{s} (w+1)  e^{- r}. \end{equation}
\item[(c)] In the case of a discrete space $E=\{x_1,\ldots, x_q\}$ we even have the more precise estimate
\begin{equation} D^{(e)}( K \left( \lambda r + (\lambda+1) \log q \right) ~|~ X,\rho_\DD,s) \leq  4 K^{1/s} w \min(1,\lambda)^{1/s}\,  e^{- r},\label{eqn:view2}\end{equation} for $r>r_0=r_0(\lambda)$.
\end{enumerate} \label{thm:e}\end{thm}

Theorem~\ref{thm:e} can be interpreted in the following way. In order to quantize with error $e^{-r}$ one needs, on average, $\lambda r$ nats to encode the jump positions, $\lambda \log N(E,\rho,e^{-r})$ nats in order to encode the increments, and another $\log N(E,\rho,e^{-r})$ nats in order to encode the initial position $X(0)$. In particular, the same result can be proved without the $\log N(E,\rho,e^{-r})$ term if the initial value of the process is deterministic.

Let us compare Theorem~\ref{thm:q} and Theorem~\ref{thm:e} in the case where $N(E,\rho,\eps)\leq q \eps^{-\gamma}$. We point out that the asserted rate of the quantization error is different to the one of the entropy coding error. As we will see below neither the quantization error bounds nor the entropy coding error bounds can be improved significantly.

Finally note that $N(E,\rho,\eps)<\infty$ for all $\eps>0$ does not necessarily imply that $w:=\sup_{x,y\in E} \rho(x,y)<\infty$ if $\rho$ does not satisfy the triangle inequality.

\smallskip
For the remainder of this subsection, let us assume that $(E,\norm{.})$ is a normed linear space with distortion measure $\rho(x,y)=\norm{x-y}$.  We assume that the jump destinations of $X$ (and thus the increments $Z^{(i)}$) are independent of the jump positions (\conda). Furthermore, assume that the increments, conditioned upon $N_X=k$, are identically distributed (not necessarily independent among each other) with the same law as the $E$-valued random variable, say, $Z^{(1)}$. Furthermore we assume that $X(0)$ is deterministic, i.e.\ that for some $x_0\in E$ $X(0)=x_0$ a.s.

\begin{thm} Under the above assumtions the following statements are true.
\begin{enumerate} \item[(a)] If
$$\gamma:=\limsup_{\eps\to 0} \frac{\log d^{(q)}(\eps\,|\,Z^{(1)},\norm{.},s)}{\log 1/\eps}  \in [0,\infty[,$$
then (\ref{eqn:q}) is valid for the newly defined $\gamma$.
 \item[(b)] If $d^{(q)}(\eps\,|\,Z^{(1)},\norm{.},s) <\infty$ for all $\eps>0$, then $$D^{(e)}\left(  K\left( \lambda r + \lambda  \log d^{(q)}(e^{-r}\,|\,Z^{(1)},\norm{.},s)\right) ~|~ X,\rho_\DD,s\right) \leq C\, e^{- r}.$$ holds with some constant $C>0$ depending on the parameters $s,K,\lambda$, and $\E||Z^{(1)}||^s$. \end{enumerate}
  \label{thm:assumed}
\end{thm}

Theorem~\ref{thm:assumed} relates the complexity of coding  $X$ to that of coding the increments.
If the assumptions of both Theorems~\ref{thm:e} and \ref{thm:assumed}
are satisfied, then the bounds of the latter theorem provide a better estimate since in general $d^{(q)}(\eps)\leq N(\eps)$ (see Lemma~\ref{lem:quantentropy}). However, note that in contrast to Theorems~\ref{thm:q} and~\ref{thm:e}, Theorem~\ref{thm:assumed} requires that the increments are identically distributed and independent of the jump positions. In case of Theorems~\ref{thm:q} and~\ref{thm:e}, this is not necessary since, by assumption, the space $E$ is sufficiently well-structured (in the sense of small metric entropy $N$).

Let us remark that the assumption in Theorem~\ref{thm:assumed} that $X(0)$ be deterministic is for simplicity only. If instead $X(0)$ is a random variable in $E$, one has to add $d^{(q)}(\eps\,|\,X(0),\rho,s)$ to the average number of nats needed to encode $X$ conditioned upon $X(0)$.

Finally we mention that one can also prove counterparts to assertions~(b) and~(c) of Theorem~\ref{thm:e} in the setting of Theorem~\ref{thm:assumed}.

\subsection{Lower bounds} \label{sec:ilower} As an illustration consider the case of a discrete space $E$, namely let $E=\{x_1, \ldots, x_q\}$, which was first studied in \cite{Vor07}. Then $N(E,\rho,\eps)\leq q$, and we thus obtain from Theorems~\ref{thm:q} and~\ref{thm:e}: \begin{equation} - \log D^{(q)}( r ~|~ X,\rho_\DD,s) \gtrsim \sqrt{\frac{2}{s}\, r \log r}\quad\text{and}\quad -\log D^{(e)}(  K r ~|~ X,\rho_\DD,s) \gtrsim \frac{r}{\lambda}.\label{eqn:dissresult}\end{equation}

Now we ask for lower bounds. Clearly, one cannot expect a non-trivial lower bound when only assuming (\ref{decrcond}). Thus, let us assume in this subsection that the jump positions constitute a Poisson point process and that \conda\ holds. In this case, we show in Theorem~\ref{thm:onoff} that the order of $D^{(q)}$ in (\ref{eqn:dissresult}) is in fact the true order on this scale. Below, in Theorem~\ref{thm:steffenmod}, we show that the order of $D^{(e)}$ is the correct one, too.

We consider a more general situation than a finite, discrete space. We only have to assume that there is sufficient uncertainty in the model in order to ensure that every jump indeed has to be encoded.

Concretely, assume that \conda\ holds and that the jump positions form a Poisson point process. Furthermore, we assume that, given the event that $k$ jumps occur ($\{N_X=k\}$), the moduli of the increments $Z_1, \ldots, Z_k$ are such that there are $\eps_0>0$ and $\delta_0>0$ (independent of $k$) such that for all $i=1,\ldots, k$, $\pr{Z_i>\eps_0\,|\,N_X=k} \geq \delta_0$. Additionally, we now impose that  $(E,\rho)$ is a {\it metric} space.

\begin{thm} Under the above assumptions, $$-\log D^{(q)}( r ~|~ X,\rho_\DD,s) \lesssim  \sqrt{\frac{2}{s}\, r \log r}.$$

In particular, for a discrete metric space $E=\{x_1, \ldots, x_q\}$, $$-\log D^{(q)}( r ~|~ X,\rho_\DD,s) \sim \sqrt{\frac{2}{s}\, r \log r}.$$ \label{thm:onoff}\end{thm}

Note that in view of (\ref{eqn:view2}) the rates for quantization error and entropy coding error must be different in case of a discrete metric space $E=\{x_1, \ldots, x_q\}$. Moreover, the order of convergence of the quantization error depends strongly on the moment $s$. In particular, one has for two distinct moments $0<s<s'$  that
$$
\lim_{r\to\infty} \frac{D^{(q)}( r ~|~ X,\rho_\DD,s)}{D^{(q)}( r ~|~ X,\rho_\DD,s')} =0.
$$
This contrasts earlier results on quantization where the same order of convergence is obtained for all moments $s>0$.


Let us consider a simple example.

\begin{exa} Let $X$ be an alternating Poisson process, i.e.\ \begin{equation} X(t)=\sum_{i=1}^{N(t)}(-1)^{i-1}, \label{eqn:alternpp}\end{equation} where $(N(t))_{t\in[0,1]}$ is a Poisson (counting) process with rate $\lambda$ (cf.\ Section~\ref{sec:cpp}) with the natural metric $|.|$ on $E=\{0,1\}$. Then Theorem~\ref{thm:onoff} yields $$-\log D^{(q)}( r ~|~ X,\rho_\DD,s) \sim \sqrt{\frac{2}{s}\, r \log r}.$$ \end{exa}

\begin{re} Recall that the assertions of the upper bounds, do not  require that $(E,\rho)$ is a metric space. The statement is valid for any distortion measure $\rho$. However, the lower bound from Theorem~\ref{thm:onoff} fails for general distortion measures. This can be seen from the following simple example. Let $E=\{ 0,1\}\cup \{ 2^{-n}, n\geq 1\}$ and $\rho(0,1)=1$, $\rho(0,2^{-n})=\rho(1,2^{-n})=2^{-n}$, $\rho(2^{-n},2^{-m})=1$ for $n\neq m$ and $n,m\geq 1$. Note that this is not a metric space.

Consider the alternating Poisson process, i.e.\ the model from (\ref{eqn:alternpp}). Then $X$ satisfies the assumptions of Theorem~\ref{thm:onoff}, in particular those for the moduli of the increments (since $Z_i=1$), but $D^{(q)}( r ~|~ X,\rho_\DD,s)=0$, for all $r\geq 0$.  \end{re}

Next, we will prove a lower bound for the entropy coding error. \begin{thm} \label{thm:steffenmod} Let $X$ be a jump process satisfying condition~(*). We assume that he jump positions $(Y_i)$ form  a Poisson point process with rate $\lambda$. Moreover, we suppose that $\rho$ defines a metric on $E$ and that a.s.\ the moduli of the jumps are bounded from below by $\eps_0>0$.  Then for all $s\geq 1$ and all sufficiently large $r$
$$ D^{(e)}( r ~|~ X,\rho_\DD,s) \geq \eps_0 C \min(1,\lambda) \, e^{-r/\lambda},$$ where $C>0$ is an absolute constant.\end{thm}

\begin{re} The lower bound in Theorem~\ref{thm:steffenmod} is actually shown to hold for the \emph{distortion rate function} $D(r~|~X,\rho,s)$ defined in Section~\ref{sec:ece}.\end{re}

We obtain the following corollary as a special case. \begin{cor} \label{cor:steffen}
Let $X$ satisfy the conditions of Theorem~\ref{thm:steffenmod}. Assume additionally that $X(0)$ is deterministic and consider the case of a discrete metric space $E=\{x_1,\ldots, x_q\}$ with $w := \max_{x,y\in E} \rho(x,y)$. Then for $s\ge 1$
$$ C_1 \eps_0 \min(1,\lambda) e^{-r/\lambda}\leq  D^{(e)}( r ~|~ X,\rho_\DD,s) \leq C_2 \,q \, w\, \min(1,\lambda)^{1/s} \, e^{-r/\lambda},$$ for large enough $r$ and absolute constants $C_1, C_2>0$.
\end{cor}

The corollary follows immediately from part (c) of Theorem~\ref{thm:e} and the remark after it and Theorem~\ref{thm:steffenmod}. This result shows that the bounds for the entropy coding error in Theorems~\ref{thm:e} and~\ref{thm:assumed} are tight.

\begin{exa} Consider again the alternating Poisson process from (\ref{eqn:alternpp}) with the natural metric $|.|$. Then Corollary~\ref{cor:steffen} yields, for all $s\geq 1$, $$ C_1 \min(1,\lambda) e^{-r/\lambda} \leq  D^{(e)}( r ~|~ X,\rho_\DD,s) \leq C_2 \min(1,\lambda)^{1/s} \, e^{-r/\lambda},$$ for large enough $r$ and absolute constants $C_1, C_2>0$. \end{exa}

\begin{exa} Let us illustrate the influence of a random initial position on our estimates. For this purpose, consider an alternating Poisson process with random initial position, i.e.\  $$ X(t)=X(0) + \sum_{i=1}^{N(t)}(-1)^{i-1+X(0)}, $$ where $X(0)$ equals $0$ and $1$ with probability $1/2$, respectively, cf.\ \cite{Vor07}. Our Theorem~\ref{thm:e}, part (c), and Theorem~\ref{thm:steffenmod} yield $$ C_1 \min(1,\lambda) e^{-r/\lambda} \leq  D^{(e)}( r ~|~ X,\rho_\DD,s) \leq C_2 \min(1,\lambda)^{1/s} 2^{1/\lambda} \, e^{-r/\lambda},$$ for all $s\geq 1$ and all large enough $r$ and absolute constants $C_1, C_2>0$. \end{exa}

\subsection{Application to compound Poisson processes in $\R^d$
} \label{sec:cpp} As an application of our results, let us determine the coding complexity of $\R^d$-valued compound Poisson processes. Recall that a L\'{e}vy process with finite L\'{e}vy measure is a compound Poisson process with the following structure, cf.\ e.g.\ \cite{bertoin}.

Let $(N(t))_{t\geq 0}$ be a Poisson (counting) process with intensity $\lambda>0$, i.e.\ let  $N(t):=\max\{ n \geq 0 \,:\,\sum_{i=1}^n e_j\leq \lambda t\}$, where $(e_j)$ are i.i.d.\ standard exponential random variables. Consider \begin{equation} X(t) = \sum_{i=1}^{N(t)} Z^{(i)},\qquad t\in[0,1[, \label{eqn:cpp} \end{equation} where the $Z^{(i)}$, $i=1,2, \ldots$, are i.i.d.\ and distributed according to any probability distribution in $\R^d$ with $\pr{Z^{(1)}=0}=0$. 
Note that this notation is consistent with the one employed above for the increments. Note furthermore that for compound Poisson processes \conda\ is satisfied.

We consider the distortion measure $$\norm{X}_1 := \int_0^1 \norm{X(t)}_\infty \, \d t,$$ which of course coincides with $\rho_\DD$ for $\rho=\norm{.}_\infty$, where as usual $\norm{x}_\infty:=\max_{i=1,\ldots, d} |x_i|$. However, one can replace $\norm{.}_\infty$ by any norm on $\R^d$, which would change only the constants.

Theorem~\ref{thm:assumed} yields the following corollary. \begin{cor} Let $X$ be a compound Poisson process as defined in (\ref{eqn:cpp}) and $s>0$. \begin{enumerate} \item[(a)]
 Assume that \begin{equation} \gamma:=\limsup_{\eps\to 0} \frac{\log d^{(q)}(\eps\,|\,Z^{(1)},\norm{.}_\infty,s)}{\log 1/\eps}  \label{eqn:dbox}\end{equation}
 is finite. Then $$- \log D^{(q)}( r ~|~ X,\rho_\DD,s) \gtrsim \sqrt{\frac{2}{s (1+\gamma)}\, r \log r}.$$ \item[(b)] Let $d^{(q)}(\eps\,|\,Z^{(1)},\rho,s) <\infty$ for all $\eps>0$. Then, for $r\geq r_0$ and a constant $C=C(s,\lambda,\E || Z^{(1)}||_\infty^s)$, we have
$$ D^{(e)}( \lambda r + (\lambda+1) d^{(q)}(e^{-r}\,|\,Z^{(1)},\rho,s) ~|~ X,\rho_\DD,s) \leq C e^{- r}. $$ \end{enumerate} \label{cor:cpp}\end{cor}

Alternatively, one can study the consequences of Theorems~\ref{thm:q} and~\ref{thm:e} if one has additional information on the range of $X$.

As for lower bounds we can apply Theorem~\ref{thm:onoff}, which gives the following.
\begin{cor} Let $X$ be a compound Poisson process as defined in (\ref{eqn:cpp}) and $s>0$. Then $$ - \log D^{(q)}( r ~|~ X,\rho_\DD,s) \lesssim  \sqrt{\frac{2}{s}\, r \log r}.$$
If additionally (\ref{eqn:dbox}) holds with $\gamma=0$, then $$ - \log D^{(q)}( r ~|~ X,\rho_\DD,s) \sim \sqrt{\frac{2}{s}\, r \log r}.$$ \end{cor}

We obtain a similar result in the case that the distribution of the increments has an absolutely continuous component. 
\begin{thm} Let $X$ be a compound Poisson process as defined in (\ref{eqn:cpp}) and $s>0$. Assume that the distribution of $Z^{(1)}$ has an absolutely continuous component. Then $$-\log D^{(q)}( r ~|~ X,\rho_\DD,s) \lesssim \sqrt{\frac{2}{s(1+d)}\, r \log r}.$$

If additionally (\ref{eqn:dbox}) holds with $\gamma=d$ then $$ - \log D^{(q)}( r ~|~ X,\rho_\DD,s) \sim  \sqrt{\frac{2}{s(1+d)}\, r \log r}.$$ \label{thm:lowerac}\end{thm}

Theorems~\ref{thm:onoff} and~\ref{thm:lowerac} show that the upper bound for the quantization rate in Theorems~\ref{thm:q} and~\ref{thm:assumed} (and thus Corollary~\ref{cor:cpp}) cannot be improved in general (for all $\gamma\in\N$).

Let us finally list a corollary of Theorem~\ref{thm:steffenmod}. \begin{cor} Let $X$ be a compound Poisson process as defined in (\ref{eqn:cpp}). Assume that $\norm{Z^{(1)}}_\infty>\eps_0$ a.s. Then, for all $s\geq 1$, $$ D^{(e)}( r ~|~ X,\rho_\DD,s) \geq \eps_0 C \min(1,\lambda) \, e^{-r/\lambda},$$ for $r>r_0$ and $C>0$ an absolute constant.
\end{cor}

The most instructive examples of the application of the results of this subsection are given now. \begin{exa} Consider a Poisson (counting) process with intensity $\lambda$, i.e.\ let $Z^{(1)}=1$. Then  $$-\log D^{(q)}( r ~|~ X,\rho_\DD,s) \sim  \sqrt{\frac{2}{s}\, r \log r}$$ and $$C_1 \min(1,\lambda) e^{-r/\lambda}\leq  D^{(e)}( r ~|~ X,\rho_\DD,s) \leq C_2 e^{-r/\lambda},$$ for $s\geq 1$, $r>r_0$,  where $C_1>0$ is an absolute constant, and $C_2>0$ depends on $s$ and $\lambda$. \label{example:PP} \end{exa}

\begin{exa} Let $Z^{(1)}$ be uniformly distributed in $[0,1]^d$. Then $$-\log D^{(q)}( r ~|~ X,\rho_\DD,s) \sim  \sqrt{\frac{2}{s(1+d)}\, r \log r}$$ and $$D^{(e)}( r ~|~ X,\rho_\DD,s) \leq C_2 e^{-r/((1+d)\lambda)},$$ for $s\geq 1$, $r>r_0$, where $C_2>0$ depends on $s$ and $\lambda$. We conjecture that the order on the right-hand side is the correct one.\end{exa}

\begin{exa} Let $Z^{(1)}$ be uniformly distributed in $C$, where $C$ is the Cantor set in $[0,1]$. Set $\gamma:=\log 2/\log 3$.  Then $$\sqrt{\frac{1}{s}\, r \log r} \gtrsim -\log D^{(q)}( r ~|~ X,\rho_\DD,s) \gtrsim  \sqrt{\frac{2}{s(1+\gamma)}\, r \log r}$$ and $$ D^{(e)}( r ~|~ X,\rho_\DD,s) \leq C_2 e^{-r/((1+\gamma)\lambda)},$$ for $s\geq 1$,s $r>r_0$, where $C_2>0$ depends on $s$ and $\lambda$. We conjecture that the orders on the right-hand side, respectively, are the correct ones.\end{exa}

The theorems and examples presented in this subsection complement results from \cite{aurzadadereich}, where general real-valued L\'{e}vy processes are studied. The main result for compound Poisson processes in that paper states that, for any compound Poisson process with $\E \log \max(|Z^{(1)}|,1) < \infty$ and all $s\geq 1$, $$\log D^{(e)} (r~|~ X,\rho_\DD,s) \approx - r.$$ No result on the quantization error for compound Poisson processes is obtained in \cite{aurzadadereich}.

Our findings also improve the results in \cite{luschgypages06}, where an upper bound for the quantization error of real-valued compound Poisson processes is obtained. In particular, it is shown that for the Poisson (counting) process and all $s\geq 1$, $$-\log D^{(q)}( r ~|~ X,\rho_\DD,s) \gtrsim \sqrt{\frac{1}{s}\, r \log r}.$$ The correct rate on this scale is given in Example~\ref{example:PP}.

\section{Upper bounds} \label{sec:pg}
In this section, we provide the proofs of the upper bounds for the quantization error and the entropy coding error stated in Theorems~\ref{thm:q},~\ref{thm:e}, and~\ref{thm:assumed}, respectively. In the proofs, the following four technical lemmas are needed.

First we prove a result on the asymptotic behaviour of a certain sum occurring in the calculations. \begin{lem} Let $c>0$. Then $$\log \left( \sum_{k=0}^\infty \frac{c^k}{k!}\, e^{-c} e^{- r / (k+1)} \right)~ \sim~ - \sqrt{2 r \log r},\qquad \text{as $r\to \infty$.}$$ \label{lem:crucialorder}\end{lem}
\begin{proof} Let $V$ be a random variable that is Poisson distributed with mean $c$. Then the term in question equals $$\log \E e^{- r / (V+1)}.$$ By the so-called de Bruijn Tauberian theorem (cf.\ \cite{bgt}, Theorem~4.12.9), considering the Laplace transform is equivalent to considering the lower tail of $(V+1)^{-1}$. Thus, consider $$\log \pr{\frac{1}{V+1} <\eps} = \log \pr{V > \frac{1}{\eps} -1 }  = \log \sum_{\{ k> \frac{1}{\eps}-1\}} \frac{c^k}{k!}\, e^{-c} \sim -\frac{1}{\eps} \log \frac{1}{\eps},$$ where we used Stirling's Formula in the last step. Using the above-mentioned Tauberian theorem returns the asserted order of the Laplace transform, including the constant.\end{proof}

Secondly, we prove a quantization result for random variables in a space $E$ with known metric entropy. This is needed in order to encode the increments of the process $X$.

\begin{lem} Let $X$ be any random variable on a space $E$. Then, for all $s>0$ and all $\eps>0$, $$D^{(q)}( \log N(E,\rho,\eps) ~|~ X,\rho,s) \leq \eps.$$ In other words, $d^{(q)}( \eps~|~ X,\rho,s) \leq N(E,\rho,\eps)$. \label{lem:quantentropy} \end{lem}
\begin{proof} For given $\eps>0$ let $\CC$ be an $\eps$-net of $(E,\rho)$. By the definition of the covering numbers, $\CC$ can be chosen to contain only $N(E,\rho,\eps)$ elements. Thus $$D^{(q)}( \log N(E,\rho,\eps) ~|~ X,d,s) \leq \left( \E \min_{a\in\CC} \rho(X,a)^s \right)^{1/s} \leq \eps.$$ \end{proof}

\begin{re} By using product quantization, it is clear that for a random variable $X$ in $E^k:=E\times \ldots \times E$ with $\rho^k(x,y) := \max_{i=1,\ldots, k} \rho(x_i,y_i)$ we have $$D^{(q)}( k \log N(E,\rho,\eps) ~|~ X,\rho^k,s) \leq \eps.$$ \label{rem:pqe}\end{re}

Essentially the same technique is applied in the proof of the next lemma. The result is comparable, but slightly more precise. This version is used to encode the jump positions. \begin{lem} Let $Y$ be any random variable in $[0,1]^k$. Then, for all $s>0$, $r\geq 0$, $$D^{(q)}( r ~|~ Y,\norm{.}_\infty,s) \leq e^{-r/k}.$$ If $Y$ is such that $Y_1\leq \ldots \leq Y_k$ almost surely then we can restrict ourselves to codebooks $\CC$ with $\hat Y_1\leq \ldots \leq \hat Y_k$ for all $\hat Y \in \CC$.\label{lem:unifde} \end{lem}

Note that this may be a fairly weak estimate in concrete cases; however, it holds for all $k\geq 1$ and all $r\geq 0$. If more is known about the distribution of $Y$, much better (asymptotic) estimates are available, cf.\ \cite{grafluschgy}, e.g.\ Theorem~6.2.

\begin{proof} Let us first consider the case $e^r=n=((m+1)/2)^k$ with $m\geq 1$. Then we can use a simple product quantizer. Namely, we set $$\CC:= \{ (y_1/m,\ldots,y_k/m)\in[0,1]^k ~:~ y_i\in\{1,3,5,\ldots\}, i=1,\ldots, k\}.$$ Then $\# \CC \leq ((m+1)/2)^k$ and thus $$D^{(q)}( \log ((m+1)/2)^k ~|~ Y,\norm{.}_\infty,s) \leq \left(\E \min_{a\in \CC} \norm{Y-a}_\infty^s\right)^{1/s} \leq m^{-1}.$$

For any $r>0$ with $e^r\geq 2^k$, there exists an $m\geq 1$ such that $(\frac{m+1}{2})^k \leq e^r < (\frac{m+2}{2})^k$. Then \begin{multline*} D^{(q)}( r ~|~ Y,\norm{.}_\infty,s) \leq  D^{(q)}( \log ((m+1)/2)^k ~|~ Y,\norm{.}_\infty,s) \\  \leq m^{-1} \leq (2 e^{r/k}-2)^{-1} \leq e^{-r/k},\end{multline*} where we used $e^r\geq 2^k$ in the last step.

Finally, for $1\leq e^r\leq 2^k$, $$ D^{(q)}( r ~|~ Y,\norm{.}_\infty,s) \leq \left(\E \norm{Y-(1/2, \ldots, 1/2)}_\infty^s\right)^{1/s}\leq 1/2 \leq e^{-r/k}.$$
\end{proof}

The last lemma can be strengthened if it is known that the random vector $Y=(Y_1, \ldots , Y_k)$ satisfies $Y_1\leq \ldots \leq Y_k$.

\begin{lem} There are absolute constants $c^*, \kappa>0$ such that, for any random variable $Y$ in $[0,1]^k$ such that almost surely $Y_1\leq \ldots \leq Y_k$, we have, for all $s>0$, \begin{equation} D^{(q)}( r ~|~ Y,\norm{.}_\infty,s) \leq \frac{\kappa}{k}\, e^{-r/k},\qquad \text{for all $r\geq c^* k$.} \label{eqn:save1} \end{equation} \label{lem:improvedlemorder} \end{lem}

\begin{proof} Let $m\geq k$ and consider $$\CC := \left\lbrace \left(\frac{y_1}{m}, \ldots, \frac{y_k}{m}\right) ~:~ 1\leq y_1\leq y_2\leq \ldots \leq y_k\leq m, y_i\in\{ 1,\ldots, m\}\right\rbrace.$$ Clearly, $$\# \CC=\binom{m+k-1}{k}.$$ Note that for any $y\in[0,1]^k$ with $y_1\leq\ldots\leq y_k$ we have $\min_{a\in\CC} \norm{y-a}_\infty \leq 1/m$. Thus, $$D^{(q)}(\log (\#\CC)~|~ Y,\norm{.}_\infty,s) \leq \frac{1}{m},$$ for any random variable $Y$ that satisfies the assumption of the lemma.

Note that, by Stirling's Formula, for some absolute constants $C_1,C_2,c^*>0$, \begin{equation} \binom{k+1+k-1}{k} =\frac{(2k)!}{(k!)^2}\leq C_1 \frac{ 2^{2k}k^{2k}e^{-2k}\sqrt{2 \pi 2 k}}{k^{2k}e^{-2k}2\pi k} \leq C_2 2^{2k}\leq e^{c^* k}. \label{eqn:lqstar}\end{equation}

Let $r\geq c^* k$. Then there is an $m\geq k$ such that  \begin{equation} \binom{m+k-1}{k} < e^r \leq \binom{m+1+k-1}{k}, \label{eqn:lqplus}\end{equation} because, as seen in (\ref{eqn:lqstar}), $$\min_{m\geq k} \binom{m+k-1}{k} = \binom{k+k-1}{k} < \binom{k+1+k-1}{k} \leq e^{c^* k} \leq e^r.$$ Thus, \begin{equation}D^{(q)}(r~|~ Y,\norm{.}_\infty,s)\leq D^{(q)}\left(\left.\log \binom{m+k-1}{k}~\right|~ Y,\norm{.}_\infty,s\right) \leq \frac{1}{m}.\label{eqn:lqtristar}\end{equation}

By (\ref{eqn:lqplus}) and Stirling's Formula, for some absolute constants $C_3,C_4>0$,
\begin{multline} e^r\leq \binom{m+1+k-1}{k}= \frac{(m+k)!}{m! k!} \\ \leq C_3 \, \frac{(m+k)^m(m+k)^k e^{-m-k} \sqrt{2\pi (m+k)}}{m^m e^{-m}\sqrt{2 \pi m}\, k^k e^{-k}\sqrt{2 \pi k}} \\ \leq C_4\, \left(1+\frac{k}{m}\right)^m \frac{(2 m)^k}{k^k} \, \sqrt{\frac{m+k}{m k}}.\label{eqn:lqdoustar}\end{multline}
Observe that $\left(1+\frac{k}{m}\right)^m\leq e^k$ and that $\frac{m+k}{m k}=\frac{1}{k}+\frac{1}{m}\leq 2$, for all $m$ and $k$. Therefore, the term in (\ref{eqn:lqdoustar}) can be estimated by $$C_5 (2 e)^k \,\frac{m^k}{k^k}\leq \kappa^k\, \frac{m^k}{k^k},$$ where $\kappa$ is an absolute constant. This implies $k e^{r/k}\leq \kappa m$ or $1/m\leq \kappa e^{-r/k}/k$. We deduce from (\ref{eqn:lqtristar}) that for any $r\geq c^* k$ (\ref{eqn:save1}) holds, as asserted.
\end{proof}

Now we can proceed with the proof of our first main result.

\begin{proof}[ of Theorem~\ref{thm:q}] Let $X_k$ be a random variable that has the distribution of $X$ conditioned upon the event that $N_X=k$, i.e.\ $X$ that has $k$ jumps. Let $Y$ be the vector in $[0,1]^k$ with the jump positions of $X_k$ (in increasing order) and $Z$ be the $E^{k+1}$-vector containing values of the process $X_k$ between the jumps (in the order corresponding to when they occur), i.e.\ the initial value and the $k$ jump destinations. Note that we can reconstruct $X_k$ completely from the vectors $Y$ and $Z$. Thus, it is sufficient to find good codebooks for $Z$ and $Y$.

Let $\delta>0$. By assumption, there is an $\eps_0=\eps_0(\delta)\in]0,1[$ such that for all $0<\eps\leq\eps_0$, \begin{equation} \log N(E,\rho,\eps) \leq (\gamma+\delta) \log 1/\eps.\label{eqn:boxasmpt} \end{equation}

Let $r\geq \log 1/\eps_0$. For $0\leq k \leq k_0:=k_0(\delta,r):=r \min(1, (\log 1/\eps_0(\delta))^{-1})-1$, let $\CC_k''$ be a codebook for $Z$ in $(E^{k+1},\rho^{k+1})$ with \begin{equation}\left( \E \min_{\hat{Z}\in \CC_k''} \rho^{k+1}(Z,\hat{Z})^s\right)^{1/s} \leq 2 e^{-r/(k+1)}. \label{eqn:gets4}\end{equation} By Remark~\ref{rem:pqe}, $\CC_k''$ can be chosen such that \begin{equation}\log \#\CC_k'' \leq (k+1) \log N(E,\rho,e^{-r/(k+1)} ) \leq  (k+1) (\gamma+\delta) \log\left( e^{r/(k+1)} \right)= (\gamma +\delta) r,\label{eqn:repl1}\end{equation} where we used (\ref{eqn:boxasmpt}) and the choice of $k_0$.

For $1\leq k \leq k_0$, let $\CC_k'$ be a codebook for $Y$ in $(\R^k,\norm{.}_\infty)$ with \begin{equation}\left( \E \min_{\hat{Y}\in \CC_k'} \norm{Y-\hat{Y}}_\infty^s\right)^{1/s} \leq 2 \left( e^{r-k}\right)^{-1/k}. \label{eqn:gets2}\end{equation} By Lemma~\ref{lem:unifde}, $\CC_k'$ can be chosen such that $\log \#\CC_k' \leq r-k$.

Define $\CC_0:=\CC_0''$. For $k\neq 0$, let $\CC_k$ be the Cartesian product of the codebooks $\CC_k'$ and $\CC_k''$. Then $\log \#\CC_k \leq r-k + (\gamma +\delta)r$ for all $0\leq k\leq k_0$.

Let us define the following notation: for any $\hat{Y}\in\CC_k'$, we set \begin{equation} F:=\bigcup_{i=1}^k \left[\hat{Y}_i,Y_i\right[ \cup \left[Y_i,\hat{Y}_i\right[ \subseteq [0,1[. \label{eqn:defnff} \end{equation} Note that on $[0,1[\setminus F$, $X$ can be reconstructed up to the error given in (\ref{eqn:gets4}). Furthermore, note that the Lebesgue measure of $F$ is less than $k\,\norm{ Y-\hat{Y}}_\infty$.

With the help of this information, we can estimate the error of approximating by $\CC_k$ when $k\neq 0$: \begin{eqnarray} \E \min_{a\in \CC_k} \rho_\DD(X_k,a)^s \notag
&=& \E \min_{a\in \CC_k} \left( \int_0^1 \rho(X_k(t),a(t)) \, \d t\right)^s \notag
\\
&\leq & C_s \E \min_{a\in \CC_k} \left( \left( \int_{F} \ldots \, \d t \right)^s + \left( \int_{[0,1[\setminus F} \ldots\, \d t \right)^s\right)\notag
\\
&\leq & C_s \E \min_{\hat{Y}\in \CC_k'} \min_{\hat{Z}\in \CC_k''} \left( \left( w k \norm{Y-\hat{Y}}_\infty \right)^s + \left( \rho^{k+1}(Z,\hat{Z}) \right)^s\right)\label{eqn:crucialchange}
\\
&=& C_s \left( (k w)^s \E \min_{\hat{Y}\in \CC_k'} \norm{Y-\hat{Y}}_\infty^s +  \E \min_{\hat{Z}\in \CC_k''} \rho^{k+1}(Z,\hat{Z})^s \right)\notag
\\
&\leq&  C_s \left( (2 k w)^s \left( e^{r-k}\right)^{-s/k} + 2^s e^{-s r/(k+1)}\right) \notag
\\
& \leq& D k^s e^{-r s /(k+1)}, \label{eqn:gets1}\end{eqnarray} having used (\ref{eqn:gets2}) and (\ref{eqn:gets4}) in the last but one step, where $D:=C_s 2^s ((e w)^s +1)$.

We define the codebook $\CC := \bigcup_{0\leq k \leq k_0} \CC_k$. Then $$\#\CC\leq \sum_{0\leq k\leq k_0} e^{r-k+(\gamma+\delta) r} \leq e^{r+(\gamma+\delta) r} \sum_{k=0}^{\infty} e^{-k} \leq e^{r+(\gamma+\delta) r +1}.$$
Thus,
\begin{multline} D^{(q)}( (1+\gamma+\delta) r +1 ~|~ X,\rho_\DD,s)^s \leq  \E \min_{a\in \CC} \rho_\DD(X,a)^s \\ \leq \sum_{0\leq k \leq k_0} \pr{N_X=k} \E \min_{a\in \CC_k} \rho_\DD(X_k,a)^s  +  \sum_{k > k_0} \pr{N_X=k} \E \min_{a\in \CC_0} \rho_\DD(X_k,a)^s.\label{eqn:reas2}
\end{multline}
Using (\ref{decrcond}), (\ref{eqn:gets1}), and the trivial fact that $\rho_\DD(X_k,a)\leq w$, the last expression is seen to be less than
\begin{eqnarray}
&&K e^{-\lambda}\left( 2^s e^{-s r} + \sum_{1\leq k\leq k_0} \frac{\lambda^k}{k!}\,  D k^s  e^{-rs/(k+1)}  +  \sum_{k> k_0} \frac{\lambda^k}{k!}\, w^s\right)
\notag \\
&=& K e^{-\lambda}\left(2^s  e^{-s r} + D \sum_{1\leq k\leq k_0} \frac{\lambda^k}{k!}\,  k^s e^{-r s / (k+1)}  +  w^s \sum_{k> k_0} \frac{\lambda^k}{k!}\, e^{-rs/(k+1)} e^{rs/(k+1)}\right)
\notag \\
 &\leq& K e^{-\lambda}\left( 2^s e^{-s r} + D \sum_{1\leq k\leq k_0} \frac{\lambda^k}{k!}\, e^{ks} e^{-r s / (k+1)}  +  \left(\frac{w}{\eps_0}\right)^s \sum_{k> k_0} \frac{\lambda^k}{k!}\, e^{-r s / (k+1)} \right)
\notag \\
&\leq& K e^{-\lambda} 2^s e^{-s r} + C_{K,s,w,\lambda,\eps_0(\delta)}\sum_{k=0}^{\infty} \frac{(e^s \lambda)^k}{k!}\, e^{-e^s \lambda}  e^{-r s /(k+1)} \label{eqn:save2}.\end{eqnarray} Recall from Lemma~\ref{lem:crucialorder} that the exponential order of the sum, when $r\to\infty$, is $$- \sqrt{2 r s \log (r s)} \sim - \sqrt{2 r s \log r}$$ and that the constant in front of it does not depend on $r$. The first term in (\ref{eqn:save2}) also has no influence. Thus, for any $\delta>0$,
$$ \limsup_{r\to \infty} \frac{\log D^{(q)}( (1+\gamma+\delta)r+1 ~|~ X,\rho_\DD,s)}{\sqrt{r \log r}} \leq - \sqrt{\frac{2}{s}}.$$ Therefore $$\limsup_{r\to \infty} \frac{\log D^{(q)}( r ~|~ X,\rho_\DD,s)}{\sqrt{r \log r}} \leq - \sqrt{\frac{2}{s(1+\gamma+\delta)}},$$ which holds for any $\delta>0$. Letting $\delta$ tend to $0$ gives the assertion.\end{proof}

\begin{proof}[ of Theorem~\ref{thm:e}] First we treat part (a). Again we condition upon the event that $k$ jumps occur. Let $X_k$ be a random variable that has the distribution of $X$ conditioned upon the event that $N_X=k$, i.e.\ that $X$ has $k$ jumps. Let, as above, $Y$ be the vector in $[0,1]^k$ with the jump positions of $X_k$ and $Z$ be the $E^{k+1}$-vector containing the values of the process $X_k$ between the jumps. Recall that one can reconstruct $X_k$ from $Y$ and $Z$, so it suffices to find good codebooks for $Y$ and $Z$.

Let $c^*,\kappa>0$ be the absolute constants from Lemma~\ref{lem:improvedlemorder}. Let $r\geq \lambda c^*$. Fix $k\geq 1$. Let $\CC_k'$ be a codebook for $Y$ in $(\R^k,\norm{.}_\infty)$  with $$\left( \E \min_{\hat{Y}\in \CC_k'} \norm{Y-\hat{Y}}_\infty^s\right)^{1/s} \leq \frac{2\kappa}{k}\, \left( e^{r k/\lambda}\right)^{-1/k} = \frac{2\kappa}{k}\, e^{-r/\lambda}.$$
By Lemma~\ref{lem:improvedlemorder} and the fact that $Y_1\leq \ldots\leq Y_k$, $\CC_k'$ can be chosen such that $\log \#\CC_k' \leq k r/\lambda$.

Furthermore, for $k\geq 0$, let $\CC_k''$ be a codebook for $Z$ in $(E^{k+1},\rho^{k+1})$ with \begin{equation}\left( \E \min_{\hat{Z}\in \CC_k''} \rho^{k+1}(Z,\hat{Z})^s\right)^{1/s} \leq 2 e^{-r/\lambda}. \label{eqn:gets4be}\end{equation} By Remark~\ref{rem:pqe}, $\CC_k''$ can be chosen such that \begin{equation} \log \#\CC_k'' \leq (k+1) \log N(E,\rho, e^{-r/\lambda} ). \label{eqn:repl2}\end{equation}

Let $\CC_k$ be the Cartesian product of the codebooks $\CC_k'$ and $\CC_k''$. Then $\log \#\CC_k\leq k r /\lambda + (k+1) \log N(E,\rho, e^{-r/\lambda} )$.

Let $F$ be defined as in (\ref{eqn:defnff}). In case $k$ jumps occur ($k\neq 0$) we approximate $X$ by a function from $\CC_k$, which gives an error of at most \begin{eqnarray} \E \min_{a\in \CC_k} \rho_\DD(X_k,a)^s \notag
&=& \E \min_{a\in \CC_k} \left( \int_0^1 \rho(X_k(t),a(t)) \, \d t\right)^s \notag
\\
&\leq & C_s \E \min_{a\in \CC_k} \left( \left( \int_{F} \ldots \, \d t \right)^s + \left( \int_{[0,1[\setminus F} \ldots\, \d t \right)^s\right)\notag
\\
&\leq & C_s \E \min_{\hat{Y}\in \CC_k'} \min_{\hat{Z}\in \CC_k''} \left( \left( w k \norm{Y-\hat{Y}}_\infty \right)^s + \left( \rho^{k+1}(Z,\hat{Z}) \right)^s\right) \label{eqn:crucialchange2}
\\
&=& C_s \left( (k w)^s \E \min_{\hat{Y}\in \CC_k'} \norm{Y-\hat{Y}}_\infty^s +  \E \min_{\hat{Z}\in \CC_k''} \rho^{k+1}(Z,\hat{Z})^s \right)\notag
\\
&\leq&  C_s 2^s \left( (w \kappa)^s e^{-rs/\lambda} + e^{-s r/\lambda}\right)\notag
\\
& \leq& C_s 2^s ((w \kappa)^s+1) e^{-r s /\lambda}. \label{eqn:lsmodif} \end{eqnarray} For $k=0$, set $\CC_0:=\CC_0''$. Then the error is less than $2e^{-rs/\lambda}$, by (\ref{eqn:gets4be}).

On the other hand, this procedure has an expected nat length of at most  \begin{multline*} K e^{-\lambda} \sum_{k=0}^\infty \frac{\lambda^k}{k!}\, \log \left( \#\CC_k \right) = K \sum_{k=0}^\infty \frac{\lambda^k}{k!}\, e^{-\lambda} \left( k r /\lambda + (k+1) \log N(E,\rho, e^{-r/\lambda} )\right) \\ = K \left( r + (\lambda+1) \log N(E,\rho, e^{-r/\lambda} )\right).\end{multline*} Therefore, similarly to (\ref{eqn:reas2}), \begin{multline}  D^{(e)}( K \left( r + (\lambda+1) \log N(E,\rho, e^{-r/\lambda} )\right) ~|~ X,\rho_\DD,s)^s \\ \leq K e^{-\lambda}\left( 2e^{-rs/\lambda} + \sum_{k=1}^\infty \frac{\lambda^k}{k!}\, C_s 2^s ((w \kappa)^s+1)  e^{-rs/\lambda} \right) \\ = K \left( 2e^{-rs/\lambda} e^{-\lambda} + C_s 2^s ((w \kappa)^s+1) e^{-rs/\lambda} \sum_{k=1}^\infty \frac{\lambda^k}{k!}\, e^{-\lambda}\right)  \leq K C_s' (w^s+1) e^{-rs/\lambda}. \label{eqn:modifdisc}\end{multline}
where $C_s'$ only depends on $s$. This yields the assertion (a).

To see (b) one only has to recall that in case the jump positions are distributed as a Poisson point process we can choose $K=1$ in (\ref{decrcond}).

Let us finally show (c). In the case of a discrete space $E=\{x_1,\ldots, x_q\}$ with $w= \max_{x,y\in E} \rho(x,y)$, we can choose $\log \# \CC_k''=(k+1)\log q$. Thus, on $[0,1[\setminus F$, no error arises. This allows to replace the right-hand side in (\ref{eqn:lsmodif}) by $(2 w \kappa)^s e^{-r s /\lambda}$. Therefore, the upper bound in (\ref{eqn:modifdisc}) becomes $K (2 w \kappa)^s (1-e^{-\lambda}) e^{-rs/\lambda}$, where $\kappa>0$ is the absolute constant from Lemma~\ref{lem:improvedlemorder}. This finishes the proof of (c). \label{rem:modifindiscretecase} \end{proof}

Note that no assumption is necessary on the correlation of the jump positions and increments.

Let us now indicate the changes that are necessary to prove Theorem~\ref{thm:assumed}.

\begin{proof}[ of Theorem~\ref{thm:assumed}] The proof carries over almost literally from Theorems~\ref{thm:q} and~\ref{thm:e}, respectively. The only differences concern the assumption on $d^{(q)}$ instead of the metric entropy $N$, the fixed initial position, and the possibly unbounded jumps.

In this case, we encode the increments instead of the jump destinations. Let $Y$ be as above, but $Z$ denote the $E^{k}$ vector with the increments, i.e.\ $Z:=(Z^{(1)},\ldots,Z^{(k)})$ with $Z^{(i)} := X(Y_i) - X(Y_i-)$. Note that we can reconstruct $X$ from $Y$ and $Z$, since we asssumed $X(0)$ to be deterministic.

The first change is to replace (\ref{eqn:repl1}) by $$\log \#\CC_k'' \leq k \log d^{(q)}(e^{-r/k} ~|~ Z^{(1)},\rho,s ) \leq  k (\gamma+\delta)\log e^{r/k} = (\gamma+\delta) r$$ in the proof for the quantization error. For the entropy coding error one has to replace (\ref{eqn:repl2}) by $$\log \#\CC_k'' \leq k \log d^{(q)}(e^{-r/\lambda} ~|~ Z^{(1)},\rho,s ).$$

The second issue concerns a certain refinement in order to deal with the possibly unbounded jumps. Here, we need that we deal with a normed space. We will show that, on average, the high jumps do not have any influence on the rate. In fact, the only modification affects (\ref{eqn:crucialchange}), where we estimate by
$$C_s \E \min_{\hat{Y}\in \CC_k'} \min_{\hat{Z}\in \CC_k''} \left[ \left(\max_{1\leq n,m\leq k}\norm{\sum_{i=1}^n Z^{(i)} - \sum_{i=1}^{m} \hat{Z}^{(i)}} k \norm{Y-\hat{Y}}_\infty \right)^s + \left( k \rho^{k}(Z,\hat{Z}) \right)^s\right],$$
which is required due to the fact that we cannot estimate by a finite diameter $w$ (modification in the first term) and the errors may add up over all the jumps, since we encode the increments and not the absolute positions (modification in the second term).

The first term can be estimated by
\begin{multline*} C_s \E \min_{\hat{Y}\in \CC_k'} \min_{\hat{Z}\in \CC_k''} \left(\sum_{i=1}^k \norm{ Z^{(i)} - \hat{Z}^{(i)}} + \sum_{i=1}^k \norm{Z^{(i)}} \right)^s k^{s} \norm{Y-\hat{Y}}_\infty^s \\ \leq C_s^2 \left(k^s \E  \min_{\hat{Z}\in \CC_k''}  \rho^{k}\left( Z, \hat{Z}\right)^s + k^{s+1} \E \norm{Z^{(1)}}^s \right) k^s \E \min_{\hat{Y}\in \CC_k'}\norm{Y-\hat{Y}}_\infty^s \\ \leq C_s^2(2^s+ \E\norm{Z^{(1)}}^s) k^{2s+1} \E \min_{\hat{Y}\in \CC_k'}\norm{Y-\hat{Y}}_\infty^s, \end{multline*}
where the last step comes from (\ref{eqn:gets4}).

This leads to an additional factor $C k^{s+1}$ in (\ref{eqn:gets1}) which has no influence on the order. Note furthermore that this argument needs that the jump positions and the increments are independent (in order to separate the expectations) and that the increments are identically distributed (as $Z^{(1)}$). It is not needed that the increments are independent among each other.

Analogously, for the proof of the entropy coding error, (\ref{eqn:crucialchange2}) is modified, which leads to an additional factor of $C k^{s+1}$ in (\ref{eqn:lsmodif}), which leaves the resulting order unchanged, but which \emph{does} change the constant.
\end{proof}

\section{Lower bound for the quantization error} \label{sec:lb}
In this section, we prove the lower bounds for the quantization error. Essentially we employ a small ball argument, i.e.\ we construct an event of not too small probability that still leaves sufficient uncertainty for the error to be large.

First we prove Theorem~\ref{thm:onoff}.

\begin{proof}[ of Theorem~\ref{thm:onoff}] Let us fix $k>0$ and $\delta>0$ (to be chosen later) and define intervals $I_j:=\left[\frac{j-1}{k}+\frac{1}{4k},\frac{j}{k}-\frac{1}{4k}\right]$, $j=1,\ldots, k$. Note that $\lambda_1(I_j)=1/(2k)$. Let $A$ be the event that $X$ has exactly $k$ jumps at $Y_1, \ldots, Y_k$, such that $Y_j \in I_j$, for all $j=1, \ldots, k$, and that the moduli of the increments are all greater than $\eps_0$. Since the $Y_i$ and $Z_i$ are independent (by \conda) and the $Y_i$ are distributed according to a Poisson point process, we have \begin{eqnarray*} \pr{A} &=&\prod_{j=1}^k \pr{\text{exactly one jump in $I_j$, $Z_j>\eps_0$}} \cdot \pr{\text{no jump in $\left[(j-1)/k,j/k\right]\setminus I_j$}} \\ &\geq&
\prod_{j=1}^k \left(\frac{\lambda}{2k} e^{-\lambda/(2k)} \delta_0 \cdot e^{-\lambda/(2k)} \right)= \left(\frac{\delta_0 \lambda}{2k}\right)^k e^{-\lambda}.\end{eqnarray*}

\noindent{\it Step 1:} Let $X_A$ be a random variable with the distribution of $X$ conditioned upon the event $A$. Then
\begin{eqnarray} D^{(q)}( r ~|~ X_A,\rho_\DD,s)^s \notag
 &=& \inf_{\log (\#\CC)\leq r} \E_{X_A} \min_{f\in\CC} \rho_\DD(X_A,f)^s \notag \\
 &\geq&  \inf_{\log (\#\CC)\leq r} \delta^s \pr{ \forall f\in\CC \,:\, \rho_\DD(X_A,f) \geq \delta } \notag \\
 &=&  \inf_{\log (\#\CC)\leq r} \delta^s \left(1 - \pr{ \exists f\in\CC \,:\, \rho_\DD(X_A,f) < \delta }\right) \notag \\
 &\geq&  \inf_{\log( \#\CC)\leq r} \delta^s \left(1 - (\#\CC) \sup_f \pr{ \rho_\DD(X_A,f) < \delta }\right) \notag \\
 &\geq&  \delta^s \left(1 - e^r \sup_{f} \pr{ \rho_\DD(X_A,f) < \delta }\right),\label{smallballe}
\end{eqnarray}
where the supremum is taken over all functions $f$ in $\DD([0,1[,E)$. For such $f$, we have \begin{multline} \pr{ \rho_\DD(X_A,f) < \delta } = \pr{ \int_0^1 \rho( X_A(t),f(t)) \, \d t < \delta } \\ \leq \pr{ \bigcap_{j=1}^k \left\lbrace\int_{I_j} \rho( X_A(t),f(t)) \, \d t < \delta \right\rbrace} = \E \,\pr{ \left. \bigcap_{j=1}^k \left\lbrace\int_{I_j} \rho( X_A(t),f(t)) \, \d t < \delta \right\rbrace  \right| Z} , \label{eqn:lacr1} \end{multline}
where $Z=(X_A(0), X_A(Y_1), \ldots, X_A(Y_k))$ is the vector with the jump destinations. By \conda, we have that, conditioned upon $Z$, the events $$\left( \left\lbrace\int_{I_j} \rho( X_A(t),f(t)) \, \d t < \delta \right\rbrace \right)_{j=1}^k$$ are independent, since each of them only depends on the jump position in the respective interval. This together with (\ref{eqn:lacr1}) shows \begin{equation} \label{eqn:lblaste} \sup_{f} \pr{ \rho_\DD(X_A,f) < \delta } \leq \sup_{f} \E \prod_{j=1}^k\pr{ \left. \int_{I_j} \rho( X_A(t),f(t)) \, \d t < \delta  \right| Z} .\end{equation}

\noindent{\it Step 2:} Now we estimate each term in the product separately. Fix $j\in\{ 1,\ldots, k\}$. Define $l_j:=\frac{j-1}{k}+\frac{1}{4k}$, i.e.\ the left end point of the interval $I_j$. Furthermore, we define $$B_j:=\left\lbrace t\in I_j : \rho(X_A(l_j),f(t)) < \eps_0/2\right\rbrace\qquad\text{and}\qquad C_j:= \left\lbrace t\in I_j : X_A(t)=X_A(l_j)\right\rbrace.$$ Then we show that \begin{equation} \int_{I_j} \rho(X_A(t),f(t))\, \d t  < \delta \qquad \Rightarrow \qquad \lambda_1( B_j \Delta C_j ) < \frac{2\delta}{\eps_0}, \label{eqn:mimic01} \end{equation} where $B_j \Delta C_j := (B_j^c \cap C_j) \cup (B_j \cap C_j^c)$. Indeed, assume that we had $\lambda_1( B_j \Delta C_j ) \geq  2\delta/\eps_0$. Then \begin{multline*} \int_{I_j} \rho(X_A(t),f(t))\, \d t \\ \geq \int_{B_j^c\cap C_j} \rho(X_A(t),f(t))\, \d t + \int_{B_j\cap C_j^c} \rho(X_A(t),f(t))\, \d t \geq \frac{\eps_0}{2}\, \lambda_1( B_j \Delta C_j ) \geq  \delta,\end{multline*} where we used the triangle inequality in the last but one step. This shows (\ref{eqn:mimic01}); and we thus have $$ \pr{\left. \int_{I_j} \rho( X_A(t),f(t)) \, \d t < \delta\right|Z} \leq  \pr{ \left. | \lambda_1(B_j) - \lambda_1(C_j) | < \frac{2\delta}{\eps_0}  \right| Z} .$$ Note that, conditioned upon $Z$, $\lambda_1(B_j)$ is a deterministic value (depending on $X_A(l_j)$ and $f$), whereas $\lambda_1(C_j)$ is a random variable that is uniformly distributed in $[0,1/(2k)]$, since the point in $I_j$ where the jump of $X_A$ occurs is uniformly distributed in $I_j$. Therefore, $$ \pr{ \left. | \lambda_1(B_j) - \lambda_1(C_j) | < \frac{2\delta}{\eps_0}  \right| Z} \leq \frac{8 \delta k}{\eps_0}.$$

\noindent{\it Step 3:}  This shows, continuing (\ref{eqn:lblaste}), that $\sup_{f} \pr{ \rho_\DD(X_A,f) < \delta } \leq (8 k \delta / \eps_0)^k$. Substituting this estimate back into (\ref{smallballe}), we obtain $$D^{(q)}( r ~|~ X_A,\rho_\DD,s)^s \geq  \delta^s \left(1 - e^r (8 \delta k / \eps_0)^k\right).$$
Therefore,
$$ D^{(q)}( r ~|~ X,\rho_\DD,s)^s \geq \pr{A} \cdot D^{(q)}( r ~|~ X_A,\rho_\DD,s)^s \geq \left(\frac{\delta_0 \lambda}{2 k}\right)^k e^{-\lambda} \delta^s \left(1 - e^r \left(\frac{8 k\delta}{\eps_0}\right)^k\right). $$
Now we can optimize $k\geq 1$ and $\delta>0$ to obtain the largest possible lower bound. We set $$\delta:=\frac{\eps_0}{8k}\, \left( \frac{1}{2}\, e^{-r} \right)^{1/k}.$$ Then the last estimate becomes
$$D^{(q)}( r ~|~ X,\rho_\DD,s) \geq \left(\frac{\delta_0 \lambda}{2k}\right)^{k/s} e^{-\lambda/s} \delta \, 2^{-1/s}.$$
We set $$k := \lfloor \sqrt{ 2 s r / \log r} \rfloor \sim \sqrt{ 2 s r / \log r}.$$ Taking logarithms of the last estimate shows that $$-\log D^{(q)}( r ~|~ X,\rho_\DD,s) \lesssim \frac{k}{s} \log k + r/k \sim
 \sqrt{ \frac{2}{s}\, r \log r},$$ as asserted. \end{proof}

The proof of Theorem~\ref{thm:lowerac} contains the same idea as the one of Theorem~\ref{thm:onoff} and carries over almost literally. Therefore, we only indicate the necessary changes.

\begin{proof}[ of Theorem~\ref{thm:lowerac}]
By assumption, $Z^{(1)}$ has an absolutely continuous component. Let $S\subseteq \R^d$ be a measurable set with $\lambda_d(S)>0$ on which $Z^{(1)}$ has a positive bounded density w.r.t.\ the Lebesgue measure and such that $0\notin S$. Define $\eps_0:=\dist(S,0)/2>0$.

This time, $A$ is defined as follows: let $A$ be the event that $X$ has exactly $k$ jumps at $Y_1, \ldots, Y_k$, such that $Y_j \in I_j$, for all $j=1, \ldots, k$, and that the corresponding increments (i.e.\ $Z^{(j)}=X(Y_j)-X(Y_j-)$) are of a height in $S$. Due to the Poissonian nature of the point process and since increments and positions are independent, we have \begin{eqnarray*} \pr{A}&=&\prod_{j=1}^k \pr{\text{exactly one jump in $I_j$}} \cdot \pr{X(Y_j)-X(Y_j-)\in S}\cdot \\
&&\qquad\qquad \cdot \pr{\text{no jump in $\left[(j-1)/k,j/k\right]\setminus I_j$}} \\ &=&
\prod_{j=1}^k \left(\frac{\lambda}{2k}  e^{-\lambda/(2k)}\cdot q_S \cdot e^{-\lambda/(2k)} \right)= \left(\frac{\lambda q_S}{2k}\right)^k e^{-\lambda},\end{eqnarray*}
where $q_S:=\pr{Z^{(1)}\in S}>0$. Regarding (\ref{smallballe}), the proof is analogous to that of Theorem~\ref{thm:onoff}. We set $Z=(X_A(Y_1), \ldots, X_A(Y_k))$ for the vector with the jump destinations. In~(\ref{eqn:lacr1}) and~(\ref{eqn:lblaste}) we estimate a bit more carefully and obtain: \begin{multline*}\pr{ \rho_\DD(X_A,f) < \delta } \\ \leq  \E \prod_{j=1}^k \pr{ \left. \int_{I_j} \rho(X_A(t),f(t)) \, \d t < \delta, \int_{(4j-1)/(4k)}^{j/k} \rho(X_A(t),f(t)) \, \d t < \delta \right|Z}.\end{multline*}

As in the proof of Theorem~\ref{thm:onoff}, the sets $B_j$ and $C_j$ are introduced and (\ref{eqn:mimic01}) is established. Let $r_j':=j/k$ and $r_j:=r_j'-1/(4k)$. Because of (\ref{eqn:mimic01}) and since $X_A(t)=X_A(r_j')=X_A(r_j)$ on $[r_j,r_j']$, the last expression is less than $$\E \prod_{j=1}^k \pr{ \left.\lambda_d(B_j\Delta C_j) < \delta, \int_{r_j}^{r_j'} \rho(X_A(r_j),f(t)) \, \d t < \delta\right| Z}. $$ Note that, conditioned upon $Z$, the events $\lambda_d(B_j\Delta C_j) < \delta$ and $\int_{r_j}^{r_j'} \rho(X_A(r_j),f(t)) \, \d t < \delta$ are independent, since the second event only depends on $Z$, i.e.\ it is deterministic. Thus the last expression equals $$\E \prod_{j=1}^k \pr{ \left. \lambda_d(B\Delta C) < \delta \right| Z} \,  \pr{\left. \int_{r_j}^{r_j'} \rho(X_A(r_j),f(t)) \, \d t < \delta\right|Z}. $$ The first term can be estimated as in the proof of Theorem~\ref{thm:onoff} by $8 \delta k/\eps_0$, which allows to estimate the last expression by
$$\left(\frac{8 \delta k}{\eps_0}\right)^k\, \E \prod_{j=1}^k \pr{\left. \int_{r_j}^{r_j'} \rho(X_A(r_j),f(t)) \, \d t < \delta\right|Z}.$$

In order to treat the second term, note that it equals \begin{multline}\pr{\int_{r_j}^{r_j'} \rho(X_A(r_j),f(t)) \, \d t < \delta,j=1,\ldots, k} \\ = \E \pr{\left.\int_{r_j}^{r_j'} \rho(X_A(r_j),f(t)) \, \d t < \delta,j=1,\ldots, k\right| Z^{(1)}, \ldots, Z^{(k-1)}}. \label{eqn:verbr}\end{multline} Note that the last condition (for $j=k$) is the only non-deterministic condition in the probability. It depends on $Z^{(k)}$, which is an $\R^d$-valued random  variable distributed as $Z^{(1)}$. By the definition of the event $A$, $Z^{(k)}$ attains values in $S$. Thus, \begin{multline*}\pr{\left.\int_{r_k}^{r_k'} \rho(X_A(r_k),f(t)) \, \d t < \delta\right| Z^{(1)}, \ldots, Z^{(k-1)}} \\ = \pr{\left.\int_{r_k}^{r_k'} \norm{ \sum_{j=1}^k Z^{(j)} - f(t)}_\infty \, \d t < \delta\right| Z^{(1)}, \ldots, Z^{(k-1)}} \\ \leq\pr{\left.\norm{ \frac{Z^{(k)}}{4k} + \int_{r_k}^{r_k'}  \sum_{j=1}^{k-1} Z^{(j)} - f(t) \, \d t }_\infty < \delta\right| Z^{(1)}, \ldots, Z^{(k-1)}},\end{multline*} where the integral is to be understood componentwise. Note that $\int_{r_k}^{r_k'}  \sum_{j=1}^{k-1} Z^{(j)} - f(t) \, \d t$ is a deterministic value in $\R^d$, conditioned upon $(Z^{(1)}, \ldots, Z^{(k-1)})$. Thus, the last term is bounded from above by $\eps_0' (8 k \delta)^d$, where $\eps_0'$ is the supremum of the density of $Z^{(k)}\deq Z^{(1)}$ in $S$. In the same way, successively the other terms can be reduced; and the expression in (\ref{eqn:verbr}) can be estimated by $\eps_0'^k (8 k \delta)^{d k}$. Therefore, $$\sup_{f} \pr{ \rho_\DD(X_A,f) < \delta } \leq ( k \delta  \eps_0'')^{k(1+d)},$$ where $\eps_0''=8\min(1/\eps_0,\eps_0')$. Continuing as in Step~3 of the proof of Theorem~\ref{thm:onoff} shows $$ D^{(q)}( r ~|~ X,\rho_\DD,s)^s \geq \pr{A} \cdot D^{(q)}( r ~|~ X_A,\rho_\DD,s)^s \geq \left(\frac{\lambda q_S}{2k}\right)^k e^{-\lambda} \delta^s \left(1 - e^r \left(k\delta \eps_0''\right)^{k(1+d)}\right). $$
This time we set $$\delta:=\frac{1}{\eps_0'' k}\, \left( \frac{1}{2}\, e^{-r} \right)^{1/(k(1+d))}.$$ Then again
$$D^{(q)}( r ~|~ X,\rho_\DD,s) \geq \left(\frac{\lambda q_S}{2k}\right)^{k/s} e^{-\lambda/s} \delta \, 2^{-1/s}, $$ where this time we set $$k := \left\lfloor \sqrt{ \frac{2 s}{1+d}\, \frac{r}{\log r}} \right\rfloor \sim \sqrt{ \frac{2 s}{1+d}\, \frac{r}{\log r}}.$$ This eventually leads to $$-\log D^{(q)}( r ~|~ X,\rho_\DD,s) \lesssim \frac{k}{s} \log k + r/(k(1+d)) \sim  \sqrt{ \frac{2}{s(1+d)}\, r \log r},$$ as asserted.
\end{proof}

\section{Lower bound for the entropy coding error} \label{sec:ece}
In this section we prove a corresponding lower bound for the entropy coding error (in fact, for the distortion rate function) for a jump process where the underlying point process is Poissonian. We use the notation from \cite{Iha93}, in particular, for the distortion rate function
$$D(r~|~X,\rho_D,s):=\inf\left\lbrace \left( \E \rho_\DD(X,\hat X)^s \right)^{1/s} : I(X;\hat X)\leq r\right\rbrace,$$
and the notion of mutual information:
$$
I(X;\hat X)= \begin{cases}
\int \log \frac{\d\P_{X,\hat X}}{\d\P_X\otimes \P_{\hat X}} \,\d\P_{X,\hat X} & \text{ if } \P_{X,\hat X}\ll \P_X\otimes \P_{\hat X}\\
\infty & \text{ otherwise.}
             \end{cases}$$
We recall that $D(r~|~X,\rho,s)\leq D^{(e)}(r~|~X,\rho,s)\leq D^{(q)}(r~|~X,\rho,s)$ for any random variable, all moments and any distortion measure. Therefore, a lower bound for $D$ immediately translates into a lower bound for $D^{(e)}$.

Let us state the assumptions of the main result of this section. We shall require that $X$ is a jump process (on the index set $[0,1[$) whose jumps form a Poisson process of intensity $\lambda>0$. Furthermore, we assume  that $\rho$ defines a metric on $E$ and that the moduli of the jumps of $X$ are a.s.\ bounded from below by a constant $\eps_0>0$.

As before, we denote by  $(Y_i)$ the jump times of the process $X$ and by $N_X$ the random number of jumps of $X$; we set $Y_0=0$. Moreover, we assume that conditioned upon $N_X=k$ the random vector $(A_i):=(X(Y_i))_{i=0}^k$ and the jump times $(Y_i)$ are independent (\conda). In the rest of this section, we prove the following stronger version of Theorem~\ref{thm:steffenmod}.

\begin{thm}\label{th0930-1}
Under the above assumptions one has
$$ D( r ~|~ X,\rho_\DD,1) \geq \eps_0 C \min(1,\lambda) \, e^{-r/\lambda},$$ where $C>0$ is an absolute constant.
\end{thm}

Let us shortly describe the idea of the proof. We relate the coding complexity of the jump process to that of the random jump times.
Controlling the complexity of the jump times by using Shannon's lower bound then leads to a lower bound in terms of a variational problem.  The proof is based on several lemmas and a particular random partition of $[0,1[$.

We denote by ${\cD_m}$ the dyadic subintervals of $[0,1[$ of the $m$-th level, that is
$$
\cD_m =\{[(j-1)2^{-m}, j2^{-m}):j=1,\dots,2^m\}
$$
We construct for any collection $t_1,\dots,t_k$ of distinct points in $[0,1[$ a finite binary tree as follows. Let $\mu=\sum_{i=1}^k \delta_{t_k}$ with $\delta_z$ denoting the  Dirac mass in $z$.
The root of the tree will be associated with the interval $[0,1[$ and it will be marked by the number $\mu([0,1[)=k$. If $k\in\{0,1\}$, then the construction ends and the root is also a leaf of the tree. If $k\geq 2$, the root of the tree is attached two children namely the two dyadic intervals of $\cD_1$ that are contained in $[0,1[$: $[0,1/2[$ and $[1/2,1[$. Again we mark each of the nodes with their corresponding masses. Each node that has mass $0$ or $1$  becomes a leaf of the tree, and for each node with mass greater than $1$ we  attach the two dyadic intervals of the next level that are contained in the interval  and we continue in analogy to above.

By the construction, each leaf contains either one or no point. We shall denote by $\pi_k(t_1,\dots,t_k):=(I_1,\dots, I_k)$ the $k$-intervals associated to the leaves with positive mass. In order to make the definition unique we arrange  the intervals in their natural order.

\begin{lem}\label{le0823-1}
Let $k\geq 1$ and  $(I_1,\dots,I_k)\in \mathrm{im}(\pi_k)$. Conditioned upon the event $\{N_X=k, \pi_k(Y_1,\dots,Y_k)=(I_1,\dots,I_k)\}$ we have that $(Y_1,\dots,Y_k) \deq (U_1,\dots,U_k)$, where $U_i$ are independent random variables that are uniformly distributed on $I_i$, respectively.
\end{lem}

\begin{proof}
First note that for any collection of distinct points $t_1,\dots,t_k\in [0,1[$ such that $\sum_{j=1}^k \ind_{I_i}(t_j)=1$ for all $i=1,\dots,k$ one retrieves
$\pi_k(t_1,\dots,t_k)=(I_1,\dots,I_k)$. On the other hand, any collection of points which yields $\sum_{j=1}^k \ind_{I_i}({t_j})\not =1$ for one $i$, induces a different tree and $\pi_k(t_1,\dots,t_k)\not=(I_1,\dots,I_k)$.

Therefore, the following two events coincide $$\{N_X=k, \pi_k(Y_1,\dots,Y_k)=(I_1,\dots,I_k)\}=\Bigl\{N_X=k, \sum_{j=1}^k\ind_{I_i}({Y_j})=1 \text{ for } i=1,\dots,k\Bigr\}.$$
Recall that the times $(Y_i)$ form a Poisson process on $[0,1[$ so that conditioned on $\{N_X=k, \pi(Y_1,\dots,Y_k)=(I_1,\dots,I_k)\}$ one has $(Y_1,\dots,Y_k)\deq (U_{1},\dots,U_{k})$, where $U_i$ are independent random variables uniformly distributed on $I_i$.
\end{proof}

\begin{lem}\label{le1002-1}
Fix $k\geq 1$, $(I_1,\dots,I_k)\in \mathrm{im}(\pi_k)$ and distinct points $a_0,\dots,a_k\in E$ with $|a_{i}-a_{i-1}|\geq \eps_0$ for $i=1,\dots,k$. Moreover, let  $\mu$ denote the distribution of a process in  $\DD([0,1[,E)$ that has jump positions at $k$ uniformly distributed times in the intervals $I_1,\dots,I_k$ and that attains the values $a_0,\dots,a_k$ in the given order.
Then
$$
D(r~|~\mu,\rho_\DD,1)\geq \eps_0 \frac{k}{2e} \Bigl(\prod_{i=1}^k |I_i|\Bigr)^{1/k} \, e^{-r/k}.
$$ \label{lem:steffen1}
\end{lem}

\begin{proof}
With slight abuse of notation we shall denote by $X=(X(t))_{t\in[0,1[}$ a $\mu$-distributed process and we let $Y_1,\dots,Y_k$ denote the ordered $k$ jump positions of $X$. Due to Lemma \ref{le0823-1} the times $Y_1,\dots,Y_k$ are independent and each $Y_i$ is uniformly distributed on $I_i$.

Now let $\hat X=(\hat X(t))_{t\in[0,1[}$ denote a $\DD([0,1[,E)$-valued reconstruction with $I(X;\hat X)\leq r$. We define $X^i_t  = a_{i-1}$ for $t<Y_i$ and $X^i_t=a_i$ for $t\geq Y_i$.
Also we set $\hat X^i_t=\hat X(t)$ for $t\in I_i$ and $\hat X^i_t=X^i_t$ otherwise. Then clearly
$$
\rho_\DD (X,\hat X)\geq \sum_{i=1}^k \rho_\DD (X^i,\hat X^i).
$$
Next, we will provide a lower bound for the right hand side in the latter inequality.

For each fixed $i=1,\dots,k$ we define $\nu_i$ to be the probability kernel of the regular conditional probability  $\P(Y_i\in\cdot  |\hat X=\cdot)$. Next we choose $\hat Y_i=\hat Y_i(\hat X)$ to be the first time $t\in[0,1[$ for which the probability $\nu_i(\hat X, [0,t])$ is greater or equal to the threshold $1/2$.

We observe that for $t\in [0,1[$
\begin{multline*}
\E[ \rho( X^i_t,\hat X^i_t)  |\hat X] \geq \pr{X^i_t=a_{i-1}|\hat X}\wedge \pr{X^i_t=a_{i}|\hat X}\left[ \rho( a_{i-1},\hat X^i_t)+\rho( a_{i},\hat X^i_t)\right]\\
 \geq \pr{ X^i_t=a_{i-1}|\hat X}\wedge \pr{ X^i_t=a_{i}|\hat X} \ \rho(a_{i-1},a_i)\\
=\eps_0\ \pr{ Y_i<t|\hat X} \wedge \pr{ Y_i\geq t|\hat X}.
\end{multline*}
Consequently,  the approximation error satisfies
\begin{align*}
\E[\rho_\DD (X^i,\hat X^i)  |\hat X]  & \geq  \eps_0  \int_0^1 \pr{ Y_i<t|\hat X} \wedge \pr{ Y_i\geq t | \hat X}   \,\d t\\
&=\eps_0 \left[ \int_0^{\hat Y_i}  \E [\ind_{\{Y_i<t\}}|\hat X] \,\d t + \int_{\hat Y_i}^1  \E [\ind_{\{Y_i\geq t\}}|\hat X] \,\d t\right]\\
&=\eps_0 \, \E\left(|Y_i-\hat Y_i| \ |\hat X\right)
\end{align*}
and one gets
\begin{align}
\E[\rho_\DD (X,\hat X) ]\geq \eps_0 \sum_{i=1}^k  \E |Y_i-\hat Y_i|. \label{eqn:abslusssatz}
\end{align}

We shall now use the Shannon lower bound to derive a lower bound for the right hand side of the latter equation. For ease of notation we write shortly $Y=(Y_1,\dots,Y_k)$ and $\hat Y=(\hat Y_1,\dots,\hat Y_k)$.
We need the notation for the continuous entropy and its conditional counterpart: for $\R^k$-valued random vectors  $Z$ and $\hat Z$ we denote
$$h(Z):= - \int \log \frac{\d \P_Z}{\d \lambda_k}\, \d \P_Z \ \text{ and } \ h(Z|\hat Z):= - \int \log \frac{\d \P_{Z|\hat Z}}{\d \lambda_k}\, \d \P_{Z,\hat Z},
$$
provided the Radon-Nikodym derivatives exist and the integrals are well-defined.

Since $\hat Y$ is $\sigma(\hat X)$-measurable we have $I(Y;\hat Y)\leq I(X;\hat X)\leq r$; so that
by the Shannon lower bound
\begin{align*}
r&\geq I(Y;\hat Y)=h(Y)-h(Y|\hat Y)\\
& = h(Y)-h(Y- \hat Y|\hat Y) \geq h(Y)-h(Y- \hat Y).
\end{align*}
In particular, $Y-\hat Y$ is absolutely continuous and its differential entropy is well-defined.
Next, we set $d:=\E ||Y-\hat Y||_{\ell_1^k}$ and estimate the term $h(Y-\hat Y)$ from above by
$$
\phi(d)= \suptwo{\text{$Z$, $\P_Z\ll \lambda_k$}}{\E\norm{Z}_{\ell_1^k}\leq d} h(Z).$$ 
Using Lemma~6.4 from \cite{dereichankirchnerimkeller} (which is based on ideas from \cite{csiszar}) one can easily show that $$\sup_{\text{$Z$, $\P_Z\ll \lambda_k$}, \sum_{i=1}^k \E |Z_i|\leq d} h(Z) = \sum_{i=1}^k \log\left( \frac{2 d e}{k}\right).$$
Consequently, $r\geq h(Y)- k \log(2ed/k)$ or, equivalently,
$$
d\geq \frac{k}{2e}  \,e^{h(Y)/k}\, e^{-r/k}.
$$
Moreover, the entropy of $Y$ satisfies
$$
h(Y)=\sum_{i=1}^k h(Y_i)= \sum_{i=1}^k \log \frac1{|I_i|}.
$$
and we conclude that
$$
d\geq \frac{k}{2e} (\prod_{i=1}^k |I_i|)^{1/k} \, e^{-r/k},
$$
which together with (\ref{eqn:abslusssatz}) shows the assertion.
\end{proof}

A crucial quantity in the latter lower bound for the distortion rate function is the length of the intervals $I_i$. Later we will use the following estimate:
\begin{lem}\label{le1002-2}
Let $t_1,\dots,t_k\in[0,1[$ denote $k$ distinct points ordered by their size and let $(I_1,\dots,I_k)=\pi_k(t_1,\dots,t_k)$. With $t_0=-\infty$ and $t_{k+1}=\infty$ we get for each $i=1,\dots,k$ that
$$
|I_i|\geq \frac12 (t_{i}-t_{i-1})\wedge (t_{i+1}-t_i).
$$\label{lem:steffen2}
\end{lem}

\begin{proof} By definition $I_i$ is the largest dyadic interval that only contains the point $t_i$ and the assertion follows since all half-open intervals of length $(t_{i}-t_{i-1})\wedge (t_{i+1}-t_i)$ that contain $t_i$ do not contain any of the other points.
\end{proof}

\begin{lem}\label{le1003-2} There exists a universal constant $\alpha_1\in\R$ and a function $\alpha_2 : \N \to \R$ such that for any $k\geq 1$ and $i\in\{1,\dots,k\}$
$$
\E[\log  [(Y_i-Y_{i-1})\wedge (Y_{i+1}-Y_i)]|N_X=k]=  \alpha_1 - \alpha_2(k).
$$
\end{lem}

\begin{proof}Let $\tilde Y_1,\dots, \tilde Y_k$ denote the order statistics of $k$ independent $[0,1[$-uniformly distributed random variables, and let $(\bar Y_i)$ denote the random jump positions of a Poisson process of intensity $1$ on $[0,\infty[$.
First let  $i\in \{1,\dots,k-1\}$
\begin{align*}
\E[\log (Y_i-Y_{i-1})\wedge (Y_{i+1}-Y_i)|N_X=k]&= \E \log \bigl((\tilde Y_i-\tilde Y_{i-1})\wedge (\tilde Y_{i+1}-\tilde Y_i)\bigr)\\
&= \E\log\bigl(  \frac{ \bar Y_i- \bar Y_{i-1}}{\bar Y_{k+1}} \wedge \frac{ \bar Y_{i+1}- \bar Y_i}{\bar Y_{k+1}}\bigr)\\
&= \E \log \bigl( ( \bar Y_1- \bar Y_{0}) \wedge (\bar  Y_{2}-\bar Y_1)\bigr)- \E[\log \bar Y_{k+1}].
\end{align*}
For the second equality see e.g.\ \cite{breiman}, Proposition~13.15.

Setting $\alpha_1:=\E \log (( \bar Y_1- \bar Y_{0}) \wedge ( \bar Y_{2}-\bar Y_1))$ and $\alpha_2(k):=\E[\log \bar Y_{k+1}]$ finishes the proof in this case.
The statement follows analogously for $i=k$.
\end{proof}

Furthermore, we will  need asymptotic estimates for
$$A:=N_X\Bigl(\prod_{i=1}^{N_X} (Y_{i}-Y_{i-1})\wedge (Y_{i+1}-Y_i)\Bigr)^{1/{N_X}}$$
and $$R_\beta:={N_X} \log_+  \beta (\prod_{i=1}^{N_X} (Y_i-Y_{i-1})\wedge (Y_{i+1}-Y_i))^{1/{N_X}}= \log_+ \prod_{i=1}^{N_X} \beta (Y_i-Y_{i-1})\wedge (Y_{i+1}-Y_i),
$$
where $\beta>0$.

\begin{lem}\label{le1003-1}
One has
$$
\E R_\beta \geq \lambda \log \beta +c
$$
for the constant $c=c(\lambda)=\lambda \alpha_1- \E[N_X \alpha_2(N_X)]\in\R$, where $\alpha_1$ and $\alpha_2$ are as in the previous lemma. Moreover,
$$
\lim_{\beta\to\infty} \beta \, \E\left(\frac {N_X}\beta \wedge A\right) = \lambda.
$$
\end{lem}

\begin{proof}
Applying Lemma \ref{le1003-2} we get
\begin{align*}
\E R_\beta &= \E[\E \bigl[R_\beta|{N_X}]]\geq \E \sum_{i=1}^{N_X} \Bigl[\log \beta +\E\bigl[\log (Y_i-Y_{i-1})\wedge (Y_{i+1}-Y_i)|{N_X}\bigr]\Bigr]\\
&= \E {N_X} \Bigl[ \log \beta +\alpha_1- \alpha_2({N_X})\Bigr]= \lambda \log \beta  +c.
\end{align*}

The second statement is an immediate consequence of the monotone convergence theorem: since $A>0$ a.s.\ one has
$$
\beta  \, \E\Bigl[ \frac {N_X}\beta \wedge A\Bigr]= \E [{N_X} \wedge \beta A]\to \E {N_X}=\lambda.
$$
\end{proof}

We are now in the position to prove Theorem~\ref{th0930-1}.

\begin{proof}[ of Theorem~\ref{th0930-1}]
Let $\hat X$ be $\DD([0,1[,E)$-valued reconstruction with $I(X;\hat X)\leq r$ for some fixed $r\geq0$.
Furthermore, we denote by \begin{multline*} G(k, (I_1,\dots, I_k),(a_0,\dots,a_k)) \\ = I(X;\hat X| {N_X}=k , \pi_k(Y)=(I_1,\dots,I_k),(A_0,\dots,A_k)=(a_0,\dots,a_k))\end{multline*} the conditional mutual information of $X$ and $\hat X$ given ${N_X}$, $\pi_{N_X}(Y)$, and $(A_0,\dots, A_{N_X})$. We consider the non-negative random variable  $R= G({N_X},\pi_{N_X}(Y),(A_0,\dots,A_k))$.

Since $({N_X},\pi_{N_X}(X),(A_0,\dots,A_k))$ is $\sigma(X)$-measurable one has
$$
r\geq I(X;\hat X)\geq I(X;\hat X| {N_X}, \pi_{N_X}(Y),(A_0,\dots,A_k))= \E R.
$$
Moreover, Lemma \ref{le1002-1} together with Lemma \ref{le1002-2} implies that
$$
\E \rho_\DD (X,\hat X)\geq  \frac{\eps_0}{4e} \E \Bigl[{N_X} \Bigl(\prod_{i=1}^{N_X}  (Y_i-Y_{i-1})\wedge (Y_{i+1}-Y_i) \Bigr)^{1/{N_X}} \, e^{-R/{N_X}}\Bigr].
$$

In order to get a lower bound for the coding error we next analyze the minimization problem
$$
\E \left[ {N_X} \left(\prod_{i=1}^{N_X}  (Y_i-Y_{i-1})\wedge (Y_{i+1}-Y_i) \right)^{1/{N_X}} \, e^{-\bar R/{N_X}}\right]=\min !
$$
where the infimum is taken over all non-negative random variables $\bar R$ satisfying $\E \bar R\leq r$.

We  let again $A={N_X} (\prod_{i=1}^{N_X}  (Y_i-Y_{i-1})\wedge (Y_{i+1}-Y_i) )^{1/{N_X}}$.
Using Lagrange multipliers one gets that for every $\beta>0$
$$
R_\beta = {N_X} \log_+ \frac{\beta A}{{N_X}}={N_X} \log_+ \beta \Bigl(\prod_{i=1}^{N_X}  (Y_i-Y_{i-1})\wedge (Y_{i+1}-Y_i) )^{1/{N_X}}\Bigr),
$$
is a minimizer when $r=r_\beta:=\E R_\beta=\E[{N_X} \log_+ \frac{\beta A}{{N_X}}]$. Moreover, elementary computations give that the corresponding minimal value in the minimization problem is
$$
d_\beta:=\E \left[A \exp \left(- \log_+ \frac{\beta A}{{N_X}}\right)\right]=\E\left(\frac {N_X} \beta\wedge A\right).
$$

For given $r\geq 0$ we now choose $\beta=\beta(r)= \exp((r-c)/\lambda)$ where $c$ is as in Lemma \ref{le1003-1}. Then
$
r=\lambda \log \beta+c\leq \E R_\beta
$
and due to the variational formula above one has
$$
D(r~|~X,\rho_\DD,1) \geq \frac{\eps_0}{4e} \E \left(\frac{{N_X}}{\beta(r)} \wedge A\right).
$$
Thus letting $r$ tend to infinity we get
$$
D(r~|~X,\rho_\DD,1) \gtrsim \frac{\eps_0}{4e}\,\frac \lambda {\beta(r)} =\frac{\eps_0}{4e} \lambda \exp(-(r-c)/\lambda).
$$
Thus, one has for all sufficiently large $r$ that
$$
D(r~|~X,\rho_\DD,1) \ge  \frac{1}{8e} \lambda e^{-c/\lambda} \eps_0 e^{-r/\lambda}= C_\lambda \eps_0 e^{-r/\lambda},
$$
where
$$C_\lambda=\frac{\lambda}{8 e}\, e^{c/\lambda}\text{ and } c=\lambda \alpha_1 - \E [ N_X \,\alpha_2(N_X)].$$
Moreover, $\alpha_1$ and $\alpha_2$ can be expressed in terms of i.i.d.\ standard exponential random variables $(e_i)$ as $\alpha_1=\E \log (e_1\wedge e_2)$ and $\alpha_2(n)=\E \log \sum_{i=1}^{n+1} e_i$, cf.\ the proof of Lemma~\ref{le1003-2}.

After some calculations (using Mathematica) one obtains $$\alpha_1 =\int_0^\infty (\log x) 2 e^{-2 x} \, \d x = -\gamma - \log 2,\qquad \alpha_2(n)=\frac{\Gamma'(n+1)}{\Gamma(n+1)},$$ where $\gamma=0.57721\ldots$ is the Euler--Mascheroni constant and $\Gamma$ is the Gamma function. Some more calculations show that $$c/\lambda =- \gamma - \log 2 - \log \lambda - \int_{\lambda}^\infty x^{-1} e^{-x}\, \d x - \frac{1-e^{-\lambda}}{\lambda}.$$ Closer analysis of this term shows that $$\lim_{\lambda\to \infty} C_\lambda = \frac{1}{8 e}\, e^{-\gamma},\qquad \lim_{\lambda\to 0} C_\lambda/\lambda = \frac{1}{8 e^2},$$ which altogether shows that $C_\lambda$ can be estimated from below by $D \min(1,\lambda)$ with some absolute constant $D>0$.
\end{proof}

\noindent {\bf Acknowledgements.} The research of Frank Aurzada was supported by the DFG Research Center \textsc{Matheon}
   ``Mathematics for key technologies'' in Berlin. Christian Vormoor was supported by the DFG Graduiertenkolleg 251.

\bibliographystyle{alpha}

\begin{thebibliography}{99}
\bibitem{dereichankirchnerimkeller}
S.\ Ankirchner, S.\ Dereich, and P.\ Imkeller,
\newblock The Shannon information of filtrations and the additional logarithmic utility of insiders,
\newblock Ann.\ Probab.\ 34 (2006) 743-778.

\bibitem{aurzadadereich}
F.\ Aurzada and S.\ Dereich,
\newblock The coding complexity of L\'{e}vy processes,
\newblock Preprint (2007), available from: http://arxiv.org/abs/0707.3040

\bibitem{bertoin}
J.\ Bertoin,
\newblock L\'evy processes, volume 121 of Cambridge Tracts in
  Mathematics,
\newblock Cambridge University Press, Cambridge, UK, 1996.

\bibitem{bgt}
N.\ H.\ Bingham, C.\ M.\ Goldie, and J.\ L.\ Teugels,
\newblock Regular Variation, volume~27 of Encyclopedia of
Mathematics and its Applications,
\newblock Cambridge University Press, Cambridge, UK, 1989.

\bibitem{breiman}
L.\ Breiman,
\newblock Probability, volume~7 of Classics in Applied Mathematics,
\newblock SIAM, Philadelphia, USA, 1993.

\bibitem{CoTho91}
T.\ M.\ Cover and  J.\ A.\ Thomas,
\newblock Elements of information theory, Wiley Series in Telecommunications, John Wiley \& Sons, Inc., New York, USA, 1991.

\bibitem{CDMR08}
J.\ Creutzig, S.\ Dereich, Th.\ M\"uller-Gronbach, and K.\ Ritter,
\newblock Infinite-dimensional quadrature and approximation of distributions,
\newblock Preprint (2007).

\bibitem{csiszar}
I.\ Csiszar,
\newblock $I$-divergence geometry of probability distributions and minimization problems,
\newblock Ann.\ Probab.\ 3 (1975) 146--158.

\bibitem{Der06a}
S.\ Dereich,
\newblock The coding complexity of diffusion processes under supremum norm distortion,
\newblock to appear in: Stochastic Process.\ Appl., available from: http://dx.doi.org/ 10.1016/j.spa.2007.07.003, 2007.

\bibitem{Der06b}
S.\ Dereich,
\newblock The coding complexity of diffusion processes under ${L}^p[0,1]$-norm distortion,
\newblock to appear in: Stochastic Process. Appl., available from: http://dx.doi.org/ 10.1016/j.spa.2007.07.002, 2007.


\bibitem{DFMS03}
S.\ Dereich, F.\ Fehringer, A.\ Matoussi and M.\ Scheutzow,
\newblock On the link between small ball probabilities and the quantization problem for Gaussian measures on Banach spaces,
\newblock J.\ Theoret.\ Probab.\ 16 (2003) 249--265.

\bibitem{dereichscheutzow}
S.\ Dereich and M.\ Scheutzow,
\newblock High resolution quantization and entropy coding for fractional Brownian motion,
\newblock Electronic J.\ Prob.\ 11 (2006) 700--722.

\bibitem{grafluschgy}
S.\ Graf and H.\ Luschgy,
\newblock Foundations of {Q}uantization for {P}robability {D}istributions, volume~1730 of Lecture Notes in Mathematics,
\newblock Springer, Berlin, Germany, 2000.

\bibitem{Iha93}
S.\ Ihara,
\newblock Information theory for continuous systems,
\newblock World Scientific, Singapore, 1993.

\bibitem{Kol68}
A.\ N.\ Kolmogorov,
\newblock Three approaches to the quantitative definition of information,
\newblock Internat.\ J.\ Comput.\ Math.\ 2 (1968) 157--168.

\bibitem{LuPa04}
H.\ Luschgy and G.\ Pag\`{e}s,
\newblock Sharp asymptotics of the functional quantization problem for {G}aussian processes,
\newblock Ann.\ Probab.\ 32 (2004) 1574--1599.

\bibitem{LuPa04b}
H.\ Luschgy and G.\ Pag\`{e}s,
\newblock Functional quantization of a class of {B}rownian diffusions: a constructive approach,
\newblock Stochastic Process.\ Appl.\ 116 (2006) 310--336.

\bibitem{luschgypages06}
H.\ Luschgy and G.\ Pag\`{e}s,
\newblock Functional quantization rate and mean pathwise regularity of processes with an application to L\'{e}vy processes,
\newblock to appear in: Ann.\ Appl.\ Probab.\ 2007.

\bibitem{PaPri05}
G.\ Pag\`es and J.\ Printems,
\newblock Functional quantization for pricing derivatives,
\newblock Universit\'e de Paris VI, LPMA no.\ 930, Preprint, 2004.

\bibitem{Vor07}
C.\ Vormoor,
\newblock High resolution coding of point processes and the Boolean model,
\newblock PhD thesis, Technische Universit\"{a}t Berlin, 2007.
\end{thebibliography}

\end{document}